\documentclass[12pt]{article}
\usepackage{amssymb,latexsym,pstricks,pst-node}

\newtheorem{theorem}{Theorem}[section]
\newtheorem{proposition}[theorem]{Proposition}
\newtheorem{corollary}[theorem]{Corollary}
\newtheorem{definition}[theorem]{Definition}

\newtheorem{lemma}[theorem]{Lemma}

\newenvironment{proof}{\paragraph{\it Proof.}}{$\square$\vskip0.4cm}
\newenvironment{remark}{\paragraph{\it Remark.}}{\vskip0.4cm}

\def\A{\mathcal{A}}
\def\Circ{{\rm Circ}(\A)}

\def\Cocirc{{\rm Circ}(\B)}

\def\B{\mathcal{B}}
\def\V{\mathcal{V}}
\def\C{\mathbb{C}}

\def\P{\mathbb{P}}
\def\H{\mathbb{H}}
\def\R{\mathbb{R}}
\def\Q{\mathbb{Q}}
\def\T{\mathbb{T}}
\def\t{\mathfrak{t}}

\def\Z{\mathbb{Z}}
\def\N{\mathbb{N}}

\newcommand{\binom}[2]{\mbox{\Large $#1 \choose #2$}}

\begin{document}
\thispagestyle{empty}

\title{Toric Hyperk\"ahler  Varieties}

\author{
 Tam\'as Hausel
\\ {\it Miller Institute for Basic Research in Science}
\\ {\it and  \ Department of Mathematics}
\\ {\it University of California at Berkeley}
\\ {\it Berkeley CA 94720, USA}
\\{\tt hausel@math.berkeley.edu}
\\
\\Bernd Sturmfels
\\ {\it Department of Mathematics}
\\ {\it University of California at Berkeley}
\\ {\it Berkeley CA 94720, USA}\\
{\tt bernd@math.berkeley.edu} }

\maketitle
\begin{abstract}
Extending work of Bielawski-Dancer \cite{BD} and Konno \cite{Ko},
we develop a theory of toric hyperk\"ahler varieties, which involves
toric geometry, matroid theory and convex polyhedra.
The framework is a detailed study of semi-projective toric varieties,
meaning GIT quotients of affine spaces by torus actions, and specifically,
of Lawrence toric varieties, meaning GIT quotients of even-dimensional
affine spaces by symplectic torus actions. A toric hyperk\"ahler variety
is a complete intersection in a Lawrence toric variety. Both varieties are
non-compact, and they share the same cohomology ring, namely, the
Stanley-Reisner ring of a matroid  modulo a linear system of parameters.
Familiar applications of toric geometry to combinatorics,
including the Hard Lefschetz Theorem and the volume  polynomials of
Khovanskii-Pukhlikov \cite{KP}, are extended to the hyperk\"ahler setting.
When the matroid is graphic, our construction gives the
toric quiver varieties, in the sense of Nakajima \cite{Na}.
\end{abstract}

\section{Introduction}

Hyperk\"ahler geometry has emerged as an important
new direction in differential and algebraic geometry, with numerous
applications to mathematical physics and representation theory.
Roughly speaking, a {\it hyperk\"ahler manifold} is a
Riemannian manifold of dimension $4n$, whose holonomy is in the unitary
symplectic group $\, Sp(n)\subset SO(4n)$. The key example
is the quaternionic  space
$\,\H^n \simeq \C^{2n} \simeq \R^{4n} $.
Our aim is to relate hyperk\"ahler geometry to
the combinatorics of convex polyhedra. We believe that this
connection is fruitful for both subjects.
Our objects of study are the {\it toric hyperk\"ahler manifolds}
of Bielawski and Dancer  \cite{BD}. They are obtained from
$\, \H^n \,$ by taking the hyperk\"ahler quotient \cite{HKLR}
by an abelian subgroup of $Sp(n)$.
Bialewski and Dancer found that the geometry and topology of
toric hyperk\"ahler manifolds is
governed by hyperplane arrangements, and
Konno \cite{Ko} gave an explicit presentation of their cohomology rings.
The present paper is self-contained and contains new
proofs for the relevant  results of \cite{BD} and \cite{Ko}.

We start out in Section \ref{git} with a discussion
of semi-projective toric varieties. This may be of independent
interest.
A toric variety $X$ is called {\em semi-projective} if
$X$  has a torus-fixed point and $X$ is projective over
its affinization $\, {\rm Spec} ( H^0(X,{\cal O}_X))$.
We show that semi-projective toric varieties are
exactly the ones which arise as GIT quotients of a complex
vector space by an abelian group.
Then we calculate the cohomology ring
of a semi-projective toric orbifold $X$.
It coincides with the cohomology of the
{\em core} of $X$, which is defined as the union
of all compact torus orbit closures. This result
and further properties of the core
are derived in Section~\ref{corevariety}.

The lead characters in the present paper are the
{\em Lawrence toric varieties}, to be introduced
in Section \ref{lawrence} as the GIT quotients of
symplectic torus actions on even-dimensional affine spaces.
They can be regarded as the ``most non-compact'' among all
semi-projective toric varieties. The combinatorics of
Lawrence toric varieties is governed by the
Lawrence construction of convex polytopes \cite[\S 6.6]{Zi}
and its intriguing interplay with matroids and
hyperplane arrangements.

In Section~\ref{konno} we define
{\em toric hyperk\"ahler varieties} as
subvarieties of Lawrence toric varieties
cut out by certain natural bilinear equations.
In the smooth case, they are shown to be biholomorphic with the
toric hyperk\"ahler manifolds of Bielawski and Dancer,
whose differential-geometric construction
is reviewed  in Section~\ref{hkquotients}
for the reader's convenience. Under this identification
the core of the toric hyperk\"ahler variety coincides
with the core of the ambient Lawrence toric variety.
We shall prove that these spaces  have
the same cohomology ring which has the following
description. All terms and symbols appearing
in Theorem \ref{main} are defined in
Sections \ref{lawrence} and \ref{konno}.

\begin{theorem} \label{main}
Let $A:\Z^n\rightarrow \Z^d$ be an epimorphism,
defining an  inclusion $\,\T^d_\R  \subset \T^n_\R  \,$  of compact tori,
and let $\,\theta  \in \Z^d\,$ be generic.
Then the following graded $\Q$-algebras are isomorphic:
\begin{enumerate}
\item the cohomology ring of the
toric hyperk\"ahler variety
$\,\,Y(A,\theta)=\H^n /\!/ \! / \! /_{(\theta,0)} \T^d_\R$,
\item the cohomology ring  of the Lawrence toric variety
$\,\, X(A^\pm,\theta)=\C^{2n}/ \! /_{\theta} \, \T^d_\R$,
\item the cohomology ring  of the
core  $\,\,C(A^\pm,\theta)$, which is the preimage of the origin under
the affinization map of either the Lawrence toric variety or the toric
hyperk\"ahler variety,
\item the quotient ring $\,\Q[x_1,\dots,x_n]/(M^*(\A) + \Circ )$, where
$M^*(\A)$ is the matroid ideal which is
 generated by squarefree monomials
representing cocircuits of $A$, and
$\Circ$ is the ideal generated by the linear forms 
that correspond to elements in the kernel of $A$.
\end{enumerate}
If the matrix $A$ is unimodular then
$X(A^\pm \!,\theta)$ and $Y(A,\theta)$
are smooth and $\Q$ can be replaced by~$\Z$.
\end{theorem}

\smallskip

Here is a simple example where all three spaces are manifolds: take
$A  :  \Z^3 \rightarrow \Z , \, (u_1,u_2,u_3) \mapsto
u_1+u_2+u_3 \,$ with $\,\theta \not = 0$.
Then $\,C(A^\pm,\theta) \,$  is the
complex projective plane $\, \P^2 $.
The Lawrence toric variety $X(A^\pm,\theta)$ is
the quotient of $\,\C^6 = \C^3 \oplus \C^3$
modulo the symplectic torus action
$\,(x,y) \mapsto (t\cdot x, t^{-1} \cdot y)$.
Geometrically, $X$ is a rank $3$ bundle over $\P^2$, visualized as an
unbounded $5$-dimensional polyhedron with a bounded $2$-face,
which is a triangle.
The toric hyperk\"ahler variety $Y(A,\theta)$  is embedded into
$X(A^\pm,\theta)$
as the hypersurface $\, x_1 y_1 +  x_2 y_2 + x_3 y_3 = 0$.
It is isomorphic to the cotangent bundle of $\P^2$.
Note that $Y(A,\theta)$ itself is not a toric variety.

For general matrices $A$,
the varieties $X(A^\pm,\theta)$ and $Y(A,\theta)$ are orbifolds, by
the genericity hypothesis on $\theta$,
and they are always non-compact.
The core  $\,C(A^\pm,\theta) \,$ is projective
but almost always reducible. Each of its irreducible
components is a projective toric orbifold.

In  Section~\ref{cogenerators} we give a dual
presentation, in terms of {\em cogenerators}, for the
cohomology ring. These cogenerators are
the volume polynomials of  Khovanskii-Pukhlikov
\cite{KP} of the bounded faces of our unbounded polyhedra.
 As an application we prove
the injectivity part of the Hard Lefschetz Theorem for
toric hyperk\"ahler varieties, which, in light of the following
corollary to Theorem \ref{main},
provides  new inequalities for the  $h$-numbers of
rationally representable matroids.

\begin{corollary}
The Betti numbers of the toric hyperk\"ahler variety
$Y(A,\theta)$ are the $h$-numbers (defined in
Stanley's book \cite[\S III.3]{St1})
of the  rank $n-d$ matroid given by the
integer matrix $A$.
\end{corollary}

The {\em quiver varieties} of Nakajima \cite{Na} are hyperk\"ahler quotients
of $\H^n$ by some subgroup $G\subset Sp(n)$ which is a product of
unitary groups indexed by a quiver (i.e.~a directed graph).
 In Section~\ref{quiver} we examine
 {\em toric quiver varieties} which arise when
 $G$ is a compact torus. They are the toric hyperk\"ahler manifolds
obtained when  $A$
is the differential $\,\Z^{\rm edges} \rightarrow \Z^{\rm vertices}\,$
of a quiver. Note that our notion of toric quiver variety
is not the same as that of Altmann and Hille \cite{AH}. Theirs
are toric and projective: in fact, they are the
irreducible components of our core $C(A^\pm,\theta)$.

We close the paper by studying two examples in detail.
First in  Section~\ref{examples} we illustrate the main
results of this paper for a
particular example of a toric quiver variety, corresponding to the
complete bipartite graph $K_{2,3}$.
In the final Section~\ref{ale} we examine the ALE spaces
of type $A_n$. Curiously, these manifolds are both
toric and hyperk\"ahler, and we show that
they and their products are the only toric hyperk\"ahler manifolds which are
toric varieties in the usual sense.

\paragraph{\bf Acknowledgment.} This paper grew out of a lecture
on toric aspects of Nakajima's quiver varieties \cite{Na} given by
the second author in the Fall 2000
Quiver Varieties seminar at UC Berkeley,
organized by the first author.
We are grateful to the participants of this seminar for their
contributions.
In particular we thank Mark Haiman, Allen Knutson and
Valerio Toledano.
 (See
{\tt www.math.berkeley.edu/$\sim$hausel/quiver/} for seminar notes.)
 We thank Roger Bielawski for drawing our attention
to Konno's work \cite{Ko}, and we thank Manoj Chari
for explaining the importance of \cite{ML} for
Betti numbers of toric quiver varieties.
Both authors were supported by
the Miller Institute for Basic Research in Science,
in the form of a Miller Research Fellowship
(1999-2002) for  the first author and a
Miller Professorship (2000-2001) for
the second author. The second author was
also supported by the National Science Foundation (DMS-9970254).

\section{Semi-projective toric varieties}
\label{git}

Projective toric varieties are associated with rational polytopes,
that is, bounded convex polyhedra with rational vertices.
This section describes toric varieties  associated with
(typically unbounded) rational polyhedra.
The resulting class of semi-projective toric varieties will be seen to
equal the GIT-quotients of affine space  $\C^n$
modulo  a subtorus of $\T_\C^n$.

Let $A=[a_1,\dots,a_n]$ be a $d \times n$-integer matrix  whose
$d \times d$-minors are relatively prime. We choose
an $n \times (n \! - \! d)$-matrix  $B=[b_1,\dots,b_n]^T$ which makes
 the following sequence exact:
\begin{eqnarray}0 \, \longrightarrow \, \Z^{n-d}
\, \stackrel{B}{\longrightarrow}  \, \Z^n \,
\stackrel{A}{\longrightarrow} \, \Z^d \,
 \longrightarrow \, 0.\label{ses}\end{eqnarray}
The choice of $B$ is equivalent to choosing a basis in $\ker(A)$.
The configuration
$\B := \{b_1,\ldots,b_n\} $ in $ \Z^{n-d}$ is said to be a {\it Gale dual}
of the given vector configuration  $\A := \{a_1,\ldots,a_n\} $ in $\Z^d$.

We denote by $\T_\C$ the complex group $\C^*$
and by $\T_\R$ the circle $U(1)$. Their Lie algebras
are denoted by $\t_\C$ and $\t_\R$ respectively.
We apply the contravariant functor ${\rm Hom}( \,\cdot \,,\T_\C)$ to
the short exact sequence (\ref{ses}). This gives a short exact sequence
of abelian groups:
\begin{eqnarray}
\label{contraC}
 1 \,\,\longleftarrow \,\, \T_\C^{n-d} \,\,\stackrel{\, B^T}{\longleftarrow}
\,\, \T_\C^n \,\, \stackrel{\, A^T}{\longleftarrow} \,\,
\T_\C^d \,\, \longleftarrow \,\, 1 .
\end{eqnarray}
Thus $\T^d_\C$ is embedded as a  $d$-dimensional
subtorus of  $ \T_\C^n$. It acts on the affine space $\C^n$.
We shall construct the quotients of this action in the sense of
{\it geometric  invariant theory} (= GIT).
The ring of polynomial functions on $\C^n$ is graded by the
semigroup $\N \A \subseteq \Z^d$:
\begin{equation}
 S \,\,\, = \,\,\, \C[ x_1,\dots, x_n ] \,
, \,\quad {\rm deg}(x_{i}) \,\, = \,\, a_i \,\in\, \N\A.
\end{equation}
A polynomial in $S$ is homogeneous if and only if it is a $\T^d_\C$-eigenvector.
For $\theta \in \N \A$, let $S_\theta$ denote the
(typically infinite-dimensional)
$\C$-vector space of homogeneous polynomials of degree $\theta$.
Note that $S_\theta$ is
a module over the subalgebra $S_0$ of degree zero polynomials
in $\,S = \bigoplus_{\theta \in \N \A} S_\theta$.
The following lemma is a standard fact in
combinatorial commutative algebra.

\begin{lemma} \label{finite}
The $\C$-algebra $S_0$ is generated by a finite set
of monomials, corresponding to the minimal generators  of the
 semigroup $\,\N^n\cap {\rm im}(B)$.
For any $\theta\in \N\A$, the graded component
$S_{\theta}$ is a finitely generated $S_0$-module,
and the ring
 $S_{(\theta)}\,= \,\bigoplus_{r=0}^\infty S_{r \theta} \,$
is a finitely generated $S_0$-algebra.
\end{lemma}

The $\C$-algebra $S_0$ coincides with the ring of invariants
$S^{\T^d_\C}$.
The $S_0$-algebra  $S_{(\theta)}$
is isomorphic to $ \bigoplus_{r=0}^\infty \, t^r \cdot
 S_{r\theta}$, and we regard it as $\N$-graded by the degree of $t$.

\begin{definition} The {\em affine GIT quotient}
of $\,\C^n\, $ by the $d$-torus $\,\T^d_\C\,$ is the affine toric variety
\begin{equation}
\label{affineGIT}
X(A,0) \,\, := \,\,
\C^n / \! /_{0} \, \T^d_\C \,\, := \,\,
{\rm Spec} \,(S^{\T^d_\C}) \,\, = \,\,
{\rm Spec} \,(S_0) \,\, = \,\,
{\rm Spec} \,\bigl(\,\C[\N^n\cap {\rm im}(B)] \,\bigr).
\end{equation}
For any  $\theta\in \N\A$,  the  {\em projective GIT quotient}
of $\,\C^n \,$ by the $d$-torus $\,\T^d_\C \,$ is  the toric variety
\begin{equation}
\label{GIT}
 X(A,\theta) \quad := \quad
\C^n / \! /_\theta \, \T^d_\C \,\, := \,\,
 {\rm Proj}\,(S_{(\theta)}) \quad =
\quad{\rm Proj}\, \bigoplus_{r=0}^\infty \, t^r \cdot S_{r\theta} .
\end{equation}
\end{definition}

Recall that the isomorphism class of any toric variety is given
by a fan in a lattice. A toric variety is a {\em toric orbifold}
if its fan is simplicial. We shall describe the fans of
the toric varieties $X(A,0)$ and $X(A,\theta)$ using the
notation  in Fulton's book \cite{Ful}.
We write $M$ for the lattice $\Z^{n-d}$ in (\ref{ses})
and  $\,N = {\rm Hom}(M,\Z) \,$ for its dual.
The torus $\,\T_\C^{n-d} \,$ in (\ref{contraC})
is identified with  $\, N \otimes \, \T_\C $. The
column vectors $\,\B = \{b_1,\ldots,b_n\} \,$
of the matrix $B^T$ form a configuration in
$\,N\,\simeq \, \Z^{n-d}$.
We write $\,{\rm pos} (\B)\,$ for the convex polyhedral cone spanned
by $\B$ in the vector space $\,N_\R\, = \, N \otimes \, \R  \,
\simeq\, \R^{n-d}$.
Note that the affine toric variety associated with  the cone
$\,{\rm pos} (\B)\, $ equals $\,X(A,0) $.

A {\em triangulation} of the configuration $\B$ is a
simplicial fan $\Sigma$ whose rays lie in
$\B$ and whose support equals $\,{\rm pos}(\B)$.
A {\em T-Cartier divisor} on $\Sigma$ is a
continuous function $\,\Psi:{\rm pos}( \B) \to \R\, $ which is
linear on each cone of $\Sigma$ and takes integer values
on $N\cap {\rm pos}(\B)$. The triangulation $\Sigma$ is
called {\em regular} if there exists a
 T-Cartier divisor $\Psi$ which is  {\em ample}, i.e.
the function $\,\Psi:{\rm pos}( \B) \to \R\, $
is convex and restricts to a different linear function
on each maximal cone of $\Sigma$. Two T-Cartier divisors $\Psi_1$ and $\Psi_2$
are {\em equivalent} if  $\,\Psi_1-\Psi_2\,$ is a linear map on
$\,{\rm pos}(\B)$, i.e. it is an element of $M$.
A {\em divisor} on $\Sigma$ is an equivalence class of T-Cartier
divisors on $\Sigma$. Since $\Psi_1$ is ample if and only if
$\Psi_2$ is ample, ampleness is well-defined for divisors $[\Psi]$.
Finally, we define a
{\em polarized triangulation} of $\B$ to be
a pair consisting of a triangulation $\Sigma$
of $\B$ and an ample divisor $[\Psi]$.

The cokernel of $\, M
\, \stackrel{B}{\longrightarrow}  \, \Z^n \,$ is identified
with $\Z^d$ in (\ref{ses}) and we call it the
{\em Picard group}. Hence $\A = \{a_1,\ldots,a_n\} $
is a vector configuration in the Picard group.
The {\em chamber complex} $\Gamma(\A)$  of $\A$
is defined to be the coarsest fan with support ${\rm pos}(\A)$
that refines all triangulations of $\A$.
Experts in toric geometry will note that
$\Gamma(\A)$ equals the {\em secondary fan} of $\B$
as in \cite{Cox-survey}.
We say that $\theta\in \N \A $
is {\em generic} if it lies in an open chamber of $\Gamma(\A)$.
Thus $\,\theta\in \N \A \,$ is generic if it is not in
any lower-dimensional cone $\,{\rm pos} \{a_{i_1},\ldots,a_{i_{d-1}} \}\,$
spanned by columns of $A$. The chamber complex $\Gamma(\A)$
parameterizes the different combinatorial types of the
convex polyhedra
$$ P_\theta \quad = \quad \bigl\{ \,u \,\in \,{\R}^n \, :\,
A u = \theta, \, u \geq 0 \bigr\} $$
as $\theta$ ranges over $\, \N \A $.
In particular, $\theta $ is generic if and only if
$P_\theta$ is $(n-d)$-dimensional
and each of its vertices has exactly
$d$ non-zero coordinates (i.e.~$P_\theta$ is simple).
A vector $\theta$ in  $\N \A$ is called an
{\em integral degree} if every vertex of
the polyhedron $\,P_\theta \,$
is a lattice point in $\Z^n$.

\begin{proposition} There is a one-to-one correspondence between
generic integral degrees $\theta$ in $\N \A$
and polarized triangulations $\bigl(\Sigma,[\Psi]\bigr)$ of $\B$.
When forgetting the polarization this correspondence
gives a bijection between open chambers of $\Gamma(\A)$
and regular triangulations $\Sigma$ of $\B$.
\label{triangulation}
\end{proposition}

\begin{proof}
Given a  generic integral degree $\theta$, we construct
the corresponding  polarized triangulation
$\,\bigl(\Sigma,[\Psi]\bigr)$. First choose any
$\psi \in \Z^n$ such that $A \psi = -\theta$.
Then consider the polyhedron
$$ Q_\psi \quad := \quad
\bigl\{ \, v \in M_{\R} \, : \,
B v \, \geq \, \psi \, \bigr\}. $$
The  map $\, v \, \mapsto \, B v - \psi \,$
is an affine-linear isomorphism from $\,Q_\psi \,$
onto $P_\theta$ which identifies the
set of lattice points $\, Q_\psi \, \cap \,M \,$
with the set of lattice points $\, P_\theta \,\cap \, \Z^n $.
The set of linear functionals which are bounded
below on $\,Q_\psi \,$ is precisely the cone
$\,{\rm pos}({\cal B}) \subset  N $. Finally, define the function
$$ \,\Psi\, : \,{\rm pos}( \B) \to \R \, , \,\,
\, w \, \mapsto \, {\rm min} \, \{ \,w \cdot v \, : \, v \in Q_\psi \, \}. $$
This is the {\em support function} of $Q_\psi$, which is piecewise-linear,
convex and continuous. It takes integer values
on $\,N \,\cap {\rm pos}(\B)\,$ because each vertex of $Q_\psi$
lies in $M$. Since $Q_\psi$ is a simple polyhedron, its
{\em normal fan} is a regular triangulation $\,\Sigma_\theta \,$ of $\B$,
and $\Psi$ restricts to a different linear function on each
maximal face of $\Sigma_\theta$. Hence
$\, \bigl(\Sigma_\theta,[\Psi]\bigr)\,$ is a polarized
triangulation of $\B$.

Conversely, if we are given a polarized triangulation
$\, \bigl(\Sigma,[\Psi]\bigr)$ of $\B$, then we define $\, \psi \, := \,
(\Psi(b_1), \ldots,\Psi(b_n)) \, \in \Z^n $,
and $\theta = - A \psi $ is the corresponding
generic integral degree in $\N \A$.
\end{proof}

\begin{theorem} \label{toricfan}
Let $\theta\in \N\A$ be a generic integral degree. Then $X(A,\theta)$ is an
orbifold and equals the toric variety $X(\Sigma_\theta)$, where
$\Sigma_\theta$ is the regular triangulation
of $\B$ given by $\theta$ as in Proposition~\ref{triangulation}.
\end{theorem}

\begin{proof}
First note that the multigraded polynomial
ring $S$ is the {\it homogeneous coordinate ring} in
the sense of Cox \cite{Cox} of the toric variety $X(\Sigma_\theta)$.
Specifically, our sequence (\ref{ses}) is precisely the second row
in (1) on page 19 of \cite{Cox}.
The irrelevant ideal $B_{\Sigma_\theta}$ of  $X(\Sigma_\theta)$
equals the radical of the ideal generated by
$\, \bigoplus_{r=1}^\infty S_{r \theta} $.
Since $\Sigma_\theta$ is a simplicial fan,
by \cite[Theorem 2.1]{Cox},
$X(\Sigma_\theta)$ is the geometric quotient
of $\C^n \backslash \V(B_{\Sigma_\theta})\,$
modulo $\, \T^d_\C$. The variety
$\,\V(B_{\Sigma_\theta})\,$ consists
of the points in  $\C^n$ which are not semi-stable
with respect to the $\T^d_\C$-action. By standard results in
Geometric Invariant Theory, the geometric quotient of
the semi-stable locus in $\C^n$ modulo $\T^d_\C$ coincides with
$\, X(A,\theta) \, = \, {\rm Proj}\,(S_{(\theta)})
\, = \, \C^n / \! /_\theta \, \T^d_\C $.
Therefore $X(A,\theta)$ is isomorphic to $X(\Sigma_\theta)$.
\end{proof}

\begin{corollary}
The distinct GIT quotients $\,X(A,\theta) \, = \, \C^n / \! /_\theta \,
\T^d_\C \,$
which are toric orbifolds are in bijection with the open chambers in
$\Gamma(\A)$, and hence with the regular triangulations of $\B$.
\end{corollary}

Recall that for  every scheme $X$ there is a canonical morphism
\begin{eqnarray}\label{pix} \,\pi_X:X\mapsto X_0 \,\end{eqnarray}
to the affine scheme $X_0={\rm Spec} (H^0(X,{\cal O}_X))$ of
regular functions on $X$.  We call a toric variety $X$ {\em semi-projective}
if  $X$  has at least one torus-fixed point and
the morphism $\pi_X$ is projective.

\begin{theorem} The following three classes of toric varieties coincide:
\label{semi}
\begin{enumerate}
\item semi-projective toric orbifolds,
\item the GIT-quotients $X(A,\theta)$
 constructed in {\rm (\ref{GIT})}
where $\theta\in \N\A$ is a  generic integral degree,
\item toric varieties $X(\Sigma)$ where
$\Sigma$ is a regular triangulation of
a set $\B$ which spans the lattice~$N$.
\end{enumerate}
\end{theorem}

\begin{proof} The equivalence of the classes 2 and 3 follows from
Theorem~\ref{toricfan}. Let $X(\Sigma)$ be a toric variety in class 3.
Since $\B$ spans the lattice, the fan $\Sigma$ has a
full-dimensional cone, and hence $X(\Sigma)$ has
a torus-fixed point. Since $\Sigma$ is simplicial,
$X(\Sigma)$ is an orbifold. The  morphism $\pi_X$
can be described as follows. The ring of global
sections $\,H^0(X(\Sigma),{\cal O}_{X(\Sigma)}) \,$
is the semigroup algebra of the semigroup
in $M$ consisting of all linear functionals on $N$
which are non-negative on the support $|\Sigma|$ of $\Sigma$.
Its spectrum is the affine toric variety
whose cone is  $|\Sigma|$.  The triangulation $\Sigma$
supports an ample T-Cartier divisor
$\Psi$.
The morphism $\pi_X$ is projective since it is induced by $\Psi$.
Hence $X(\Sigma)$ is in class 1.

Finally, let $X$ be any semi-projective toric orbifold.
It is represented by a fan $\Sigma$ in a lattice $N$.
The fan $\Sigma$ is simplicial since $X$ is an orbifold, and
$|\Sigma|$ spans $N_\R$ since $X$ has at least one fixed point.
Since the morphism $\pi_X$ is projective, the fan
$\Sigma$ is a regular triangulation of a subset
$\B'$ of $|\Sigma|$ which includes the rays of $\Sigma$.
The set $\B'$ need not span the lattice $N$. We choose
any superset $\B$ of $\B'$ which is contained in
$\,{\rm pos}(\B') = |\Sigma|\,$ and which spans the lattice $N$.
Then $\Sigma$ can also be regarded as a regular
triangulation of $\B$, and we conclude that
$X$ is in class 3.
\end{proof}

\begin{remark}
1. The passage from $\B'$ to $\B$ in the last step means that
any GIT quotient of $\C^{n^\prime}$ modulo any abelian
subgroup of $\T_\C^{n^\prime}$ can be rewritten as a
GIT quotient of some bigger affine space $\C^n$ modulo
a subtorus of $\T_\C^{n}$. This construction applies
in particular when the given abelian group is finite,
in which case the initial subset $\B'$  of $N$
is linearly independent.

2. Our proof can be  extended to show the following:
if $X$ is any toric variety where the morphism $\pi_X$ is projective
then $X$ is the product of  a semi-projective toric variety
and  a torus.

3. The affinization map  (\ref{pix}) for $X(A,\theta)$ is by definition
the canonical map to $X(A,0)$.
\end{remark}

A triangulation $\Sigma$ of a subset $\B$ of $N \simeq \Z^{n-d}$
is called {\em unimodular} if every maximal cone of $\Sigma$
is spanned by a basis of $N$. This property holds if and only if
$X(\Sigma)$ is  a toric manifold (= smooth toric variety).
We say that a vector $\theta$ in $\N \A$ is
a {\em smooth degree} if
$\, C^{-1} \cdot \theta \, \geq \, 0 \,$ implies $\,
{\rm det}(C) = \pm 1 \,$
for every non-singular $d \times d$-submatrix $C$ of $A$. Equivalently,
the edges at any vertex of the 
polyhedron $P_\theta$ generate ${\rm ker}_\Z A\cong\Z^{n-d}$.
From Theorem \ref{semi} we conclude:

\begin{corollary} The following three classes of smooth
toric varieties coincide:
\begin{enumerate}
\item semi-projective toric manifolds,
\item the GIT-quotients $X(A,\theta)$
 constructed in {\rm (\ref{GIT})}
where $\theta\in \N\A$ is a  generic smooth degree,
\item toric varieties $X(\Sigma)$ where
$\Sigma$ is a regular unimodular triangulation of
a spanning set $\B \subset N$.
\end{enumerate}
\label{semismooth}
\end{corollary}

\begin{definition} \label{uni}
The matrix $A$ is called {\em unimodular} if  the following
equivalent conditions hold:
\begin{itemize}
\item all non-zero
$d \times d$-minors of $A$ have the same absolute value,
\item all $\,(n \! - \!d) \times (n \! - \! d)$-minors of
the matrix $B$ in (\ref{ses}) are $-1$, $0$ or $+ 1$,
\item every triangulation of $\B$ is unimodular,
\item every vector $\theta $ in $\N \A$ is an integral degree,
\item every vector $\theta $ in $\N \A$ is a smooth degree.
\end{itemize}
\end{definition}

\begin{corollary} \label{key}
For $A$ unimodular, every GIT quotient
$X(A,\theta)$ is a semi-projective toric manifold, and
the distinct smooth quotients $ X(A,\theta) $
are in bijection with the open chambers in $\Gamma(\A)$.
\end{corollary}

Every affine toric variety has a natural moment map
onto a polyhedral cone, and every projective toric
variety has a moment map onto a polytope. These
are described in Section 4.2 of \cite{Ful}. It
is straightforward  to extend this description
to semi-projective toric varieties.  Suppose
that the $S_0$-algebra $S_{(\theta)}$ in Lemma
\ref{finite} is generated by  a set of
$m+1$ monomials in $S_\theta$, possibly after
replacing $\theta$ by a multiple in the
non-unimodular case. Let  $\P_\C^m$ be the projective
space whose coordinates are  these monomials.
Then, by definition of ``Proj'',  the toric variety $X(A,\theta)$
is embedded as a closed subscheme
in the product  $\, \P_\C^m \times  {\rm Spec}(S_0)  $.
We have an action of
the $(n-d)$-torus  $\,\T^n_\C/\T^d_\C\,$ on $\P_\C^m$,
since $S_{\theta}$ is an eigenspace of $\T^d_\C$.
This  gives rise to a moment map
$\,\mu_1 \, : \,\P_\C^m \rightarrow \R^{n-d}$, whose image
is a convex polytope. Likewise, we have the affine moment map
$\, \mu_2 \, : \, {\rm Spec}(S_0) \rightarrow \R^{n-d}\,$ whose
image is the cone polar to  $\,{\rm pos} (\B)$.
This defines the moment map
\begin{equation}
\label{semimoment} \mu \, : \, X(A,\theta) \,\subset\,
 \P_\C^m \times  {\rm Spec}(S_0)  \, \rightarrow \, \R^{n-d},
\,\,\,
(u,v) \mapsto \mu_1(u) + \mu_2(v).
\end{equation}
The image of $\,X(A,\theta)\,$ under the moment map $\mu$
is the polyhedron $P_{\theta} \,\simeq \,Q_\psi $, since
the convex hull of its vertices equals the image of $\mu_1$
and the cone $P_0 \simeq Q_0$  equals the image of $\mu_2$.

Given an arbitrary fan  $\Sigma$ in $N$,
Section 2.3 in \cite{Ful} describes
how a one-parameter subgroup $\lambda_v$,
given by $v \in N$, acts on  the toric variety $X(\Sigma)$.
Consider any point
$x$ in $X(\Sigma)$ and let $\gamma \in \Sigma$ be the
unique cone such that $x$ lies in the orbit $O_\gamma$. The orbit $O_\gamma$
is fixed by the one-parameter subgroup $\lambda_v$ if and only
if $v$ lies in the $\R$-linear span $\R \gamma$ of $\gamma$.
Thus the irreducible
components $F_i$ of the fixed point locus
of the $\lambda_v$-action on $X(\Sigma)$  are the
orbit closures $\overline{O}_{\sigma_i}$ where
$\sigma_i$ runs over all cones in $\Sigma$ which are minimal
with respect to the property $\,v \in \R \sigma_i$.

The closure of $O_\gamma$ in $X(\Sigma)$ is
the toric variety $X({\rm Star}(\gamma))$ given by
the quotient fan  $\,{\rm Star}(\gamma) \,$ in
$\,N(\gamma) = N / (N \cap \R \gamma )$; see \cite[page 52]{Ful}.
From this we can  derive the following lemma.

\begin{lemma}
For $v\in N$ and $x \in O_\gamma$
the limit $\,\lim_{z \rightarrow 0} \lambda_v (z) \,x \, $
exists and lies in
 $F_i=\overline{O}_{\sigma_i}$ if and only if $\gamma \subseteq \sigma_i$ is a
face and the image of $v$ in
$N_\R/{\R \gamma}$ is in the relative interior of  $\sigma_i/{\R \gamma}$.
\label{limit}\end{lemma}

The set of all faces $\gamma$ of $\sigma_i$ with this
property is closed under  taking intersections and hence this set
has a unique minimal element.
We denote this minimal element by $\tau_i$.
Thus if we denote
$$ U^v_i
\quad = \quad
\bigl\{ \, x \in X(\Sigma) \, : \,
\lim_{z \rightarrow 0} \lambda_v (z) \,x \, \hbox{
exists and lies in } F_i \, \bigr\}, $$
or just $U_i$ for short, then this set
decomposes as a union of orbits as follows:
\begin{equation}
\label{describeui}   U_i
\quad = \quad \cup_{\tau_i\subseteq \gamma \subseteq \sigma_i}
O_\gamma.
\end{equation}
In what follows we further suppose  $v\in |\Sigma|$. Then Lemma~\ref{limit} implies
$X(\Sigma)=\cup_i U_i$, which is  the {\em Bialynicki-Birula
decomposition} \cite{BB2} of
the toric variety with respect to the
one-parameter subgroup $\lambda_v$.

We now apply this to our semi-projective
toric variety $X(A,\theta)$ with fan $\Sigma = \Sigma_\theta$.
The moment map $\mu_v$ for the circle action induced by 
$\lambda_v$ is given by the inner product
 $\mu_v(x)=\langle v ,\mu(x)\rangle$
with $\mu$  as in  (\ref{semimoment}).
We relabel the fixed components $F_i$
according to the values of this moment map, so that
\begin{eqnarray}\label{order} \mu_v(F_i)<\mu_v(F_j) \mbox{ implies }
 i<j.\end{eqnarray}
Given this labeling, the distinguished faces $\tau_i \subseteq \sigma_i$
have the following important property:
\begin{eqnarray} \label{propertystar} \tau_i \subseteq \sigma_j
\,\, \mbox{ implies } \,\, i\leq j.  \end{eqnarray}
This generalizes the property $(*)$
in \cite[Chapter 5.2]{Ful}, and it is  equivalent to
 \begin{eqnarray}
U_j \,\, \mbox{ is closed in }
\,\, U_{\leq j}=\cup_{i\leq j} U_i. \label{crucial}
\end{eqnarray}
This means that the Bialynicki-Birula decomposition
of $X(A,\theta)$ is {\em filtrable} in the sense of \cite{BB2}. Now we are able to prove the following
result, which is well-known in the projective case.

\begin{proposition}
\label{intcoho}
The integral cohomology of a smooth semi-projective  toric
variety $X(A,\theta)$ equals
$$H^*(X(A,\theta);\Z) \quad \cong
\quad \Z [x_1,x_2,\dots, x_n]/({\rm Circ}(\A)+I_\theta), $$
where $I_{\theta}$ is the Stanley-Reisner ideal of the
simplicial fan $\Sigma_\theta$, i.e. $I_\theta$ is generated by square-free
monomials $\,x_{i_1} x_{i_2} \cdots x_{i_k}\,$
 corresponding to non-faces 
$\{b_{i_1},b_{i_2},{\dots},b_{i_k}\}$ of $\Sigma_\theta$,
and ${\rm Circ}(\A)$ is the {\em circuit ideal}
$${\rm Circ}(\A) \quad := \quad
\langle \,\sum _{i=1}^n \lambda_i x_i  \,\,| \,\,\lambda \in \Z^n,
A \cdot \lambda =0 \,\rangle. $$
\end{proposition}

\begin{proof}
Let $D_1,D_2,\ldots, D_n$ denote the divisors corresponding to the
rays $b_1,b_2,\ldots,b_n$ in $\Sigma_\theta$. The cohomology class of
any torus orbit closure $\overline{O}_\sigma$ can be expressed in terms of the
$D_i$'s, namely if the rays in $\sigma$ are
$b_{i_1},b_{i_2},\dots,b_{i_k}$, then
 $\,[\overline{O}_\sigma] = [D_{i_1}][D_{i_2}]\cdots [D_{i_k}]$.
Following the reasoning in \cite[Section 5.2]{Ful},
we first prove that certain torus orbit closures linearly span
$H^*(X(A,\theta);\Z)$ and hence the cohomology
classes $[D_1],[D_2],\dots, [D_n]$ generate $H^*(X(A,\theta);\Z)$
as a $\Z$-algebra.

We choose  $v\in |\Sigma|$ to be generic, so that each $\sigma_i$ is
$(n-d)$-dimensional and each $F_i$ is just a point.
Then (\ref{describeui})
shows that $U_i$ is isomorphic with the affine space $\C^{n-k_i}$, where $k_i={\rm dim}(\tau_i)$.

We set $U_{\leq j}=\cup_{i\leq j} U_i$ and $U_{<j}=\cup_{i< j} U_i$.
Note that $U_j$ is closed in  $U_{\leq j}$. Thus
writing down the cohomology long exact sequence of the pair
$(U_{\leq j},U_{<j})$, we can
show by induction on $j$ that
the cohomology classes of the closures of the cells $U_i$
generate $H^*(X(A,\theta);\Z)$ additively. Because the closure of a cell
$U_i$ is the closure of a torus orbit, it follows that
the cohomology classes $[D_1],[D_2],\dots, [D_n]$
generate $H^*(X(A,\theta);\Z)$. Thus sending $x_i \mapsto [D_i]$
defines a surjective ring map
$\, \Z[x_1,\dots,x_n] \rightarrow H^*(X(A,\theta);\Z) $,
whose kernel is seen to contain $\,{\rm Circ}(\A)+I_\theta$.
That this is precisely the kernel follows from the ``algebraic moving lemma''
of \cite[page 107]{Ful}.
\end{proof}

A similar proof  works
with $\Q$-coefficients when $X(A,\theta)$ is not smooth
but just an orbifold.

\begin{corollary}
\label{Qcoho}
The rational cohomology ring of a semi-projective  toric
orbifold $X(A,\theta)$ equals
$$H^*(X(A,\theta);\Q) \quad \cong
\quad \Q [u_1,u_2,\dots, u_n]/({\rm Circ}(\A)+I_\theta). $$
\end{corollary}

In light of Corollary~\ref{Qcoho}, the Betti numbers of
$X(A,\theta)$ satisfy $b_{2i}=h_i(\Sigma_\theta)$, where $h_i(\Sigma_\theta)$
are the $h$-numbers of the Stanley-Reisner ideal $I_\theta$, cf. \cite[Section III.3]{St1}.
This observation
leads to the following result.

\begin{corollary}
\label{semibetti}
 If $f_i(P^{bd}_\theta)$ denotes the number of $i$-dimensional
bounded faces of $P_\theta$ then the Betti numbers of the
semi-projective toric orbifold $X(A,\theta)$ are
given by the following formula: \begin{eqnarray}
\label{alternatingsum}
b_{2k}
\quad = \quad
{\rm dim}_\Q H^{2k}(X(A,\theta);\Q)
\quad = \quad
\sum_{i=k}^{n-d}(-1)^{i-k}
\binom{i}{k} f_i(P^{bd}_\theta).  \end{eqnarray}
\end{corollary}

\begin{proof}
Lemma 2.3 of
\cite{St3} implies
that
\begin{equation}
\label{stanleysum}
\sum_{i=0}^{n-d} h_i(\Sigma_\theta) \cdot x^i
\,\,=\,\, \sum_{\sigma\in \Sigma_\theta
\backslash \partial \Sigma_\theta} (x-1)^{n-d-{\rm dim}(\sigma)},
\end{equation}
where $\partial \Sigma_\theta$ denotes the boundary
 of $\Sigma_\theta$.
Hence the right hand sum is over all interior cones $\sigma$ of
the fan
$\Sigma_\theta$. These cones are in
order-reversing bijection with
the  bounded faces of $P_\theta$.
Hence  (\ref{stanleysum}) is the sum of
$\, (x-1)^{{\rm dim}(F)} \,$ where $F$ runs over all bounded faces of
$P_\theta$. This proves
(\ref{alternatingsum}).
\end{proof}

\section{The core of a toric variety}
\label{corevariety}

The proof of Corollary~\ref{semibetti} shows the importance of interior cones
of $\Sigma_\theta$.  They are the ones for which the
closure of the corresponding torus orbit in $X(A,\theta)$
is compact.  This suggests the following

\begin{definition} \label{core}
The {\em core} of a semi-projective
toric variety $X(A,\theta)$ is
$C(A,\theta)
=\cup_{\sigma\in \Sigma_\theta \backslash \partial \Sigma_\theta}O_\sigma$.
Thus the core $\,C(A,\theta)\,$ is the union of all compact
torus orbit closures in $\,X(A,\theta)$.
\end{definition}

\begin{theorem} The core of a semi-projective toric orbifold
$X(A,\theta)$ is the inverse image of the origin under the
canonical projective morphism $X(A,\theta) \rightarrow X(A,0)$
as in (\ref{pix}).
It also equals the inverse image of the bounded faces of
the polyhedron $P_\theta$ under the moment map
(\ref{semimoment}) from
$X(A,\theta) $ onto $P_\theta$. In particular,
the core of $X(A,\theta)$ is a union of
projective toric orbifolds.
\label{coredescr}
\end{theorem}

\begin{proof}
On the level of fans,
the  toric morphism $X(A,\theta)  \rightarrow X(A,0)$ corresponds to
 forgetting the triangulation of the
cone $\, |\Sigma|  = {\rm pos}(\B)$.
It follows from the description of toric morphisms in Section 1.4 of \cite{Ful}
that the inverse image of the origin is the union of
the orbit closures corresponding to interior faces of
$\Sigma$. This was our first assertion.
Each face of a simple polyhedron is a simple polyhedron,
and each bounded face is a simple polytope. If $\sigma$
is the interior cone of $\Sigma$ dual to a bounded face of $P_\theta$
then the corresponding orbit closure is the
projective toric orbifold $X({\rm Star}(\sigma))$. The core $C(A,\theta)$
is the union of these  orbifolds.
\end{proof}

We fix a generic vector $v\in {\rm int}|\Sigma|$.
Then the $F_i$ above are  points  and lie in $C(A,\theta)$.
 In what follows
we shall study the action
of the one-parameter subgroup $\lambda_v$ on the core
$C(A,\theta)$. We define
$$D_i \, = \, U^{-v}_i \quad = \quad
\bigl\{ \, x \in X(A,\theta) \, : \,
\lim_{z \rightarrow \infty} \lambda_v (z) \,x \, \hbox{
exists and equals } F_i \, \bigr\}. $$
Lemma~\ref{limit} implies that this gives a
decomposition of the core: $C(A,\theta)=\cup_i D_i$. The
closure $\overline{D}_i$ is a projective toric orbifold, and it is the preimage
of a bounded face of $P_\theta$ via the moment map
(\ref{semimoment}).
If we now introduce an
ordering as in (\ref{order}) then the counterpart of (\ref{crucial})
is the following: \begin{eqnarray} \label{crucialdi}
D_{\leq j}\,\, = \,\,\cup_{i\leq j} D_i \,\,\mbox{ is compact}.\end{eqnarray}
This property of the decomposition  $C(A,\theta)=\cup_i D_i$
translates into a non-trivial statement about the
convex polyhedron $P_\theta$. Let  $P^{bd}_\theta$
denote the {\em bounded complex}, that is, the polyhedral
complex consisting of all bounded faces of $P_\theta$.
Let $P_j$ denote the bounded face of $P_\theta$
corresponding to $\overline{D}_j$,
and let $p_j$ denote the vertex of
$P_\theta$ corresponding to $F_j$.
Then $\,P_{\leq j} \, = \,\cup_{i\leq j} P_i \,$ is a subcomplex
of the bounded complex $P^{bd}_\theta$, and
$\, P_{\leq j} \backslash P_{< j}\,$ consists precisely
of  those faces of $P_j$ which contain $p_j$.
This property is called {\em star-collapsibility}.
It implies that $P_{< j}$ is a deformation retract of
$ P_{\leq j}$
and in turn that
$P^{bd}_\theta$ is contractible.  The contractibility
also follows from \cite[Exercise 4.27 (a)]{BLSWZ}. In summary we have proven
the following result.

\begin{theorem}
\label{starcollapsible}
The bounded complex $P^{bd}_\theta$ of $P_\theta$ is
star-collapsible; in particular, it
is contractible.
\end{theorem}

This theorem implies that the core of any semi-projective
toric variety is connected, since
$C(A,\theta)$ is the  preimage of the bounded complex
$P^{bd}_\theta$ under the continuous moment map.
Moreover, since the cohomology of $P^{bd}_\theta$ vanishes,
the bounded complex  does not contribute
to the cohomology of $C(A,\theta)$.
This fact is expressed in the following proposition,
which will be crucial in Section~\ref{cogenerators}.

\begin{proposition}
\label{littlemma}
 Let $C(A,\theta)$ be
the core of a semi-projective toric orbifold
and consider a class $\alpha$ in $ \, H^*(C(A,\theta);\Q)$. If
$\alpha$ vanishes on every irreducible component
of $C(A,\theta)$ then $\alpha = 0$.
\end{proposition}

\begin{proof}
Let $v\in {\rm int}|\Sigma|$, $F_i$ and $D_i$ as above.
We prove by induction on $j$ that
\begin{eqnarray} \mbox{ if }
\alpha\in H^*(D_{\leq j};\Q) \mbox{ and } \alpha\mid_{\overline{D}_i}
\, = \,0
\mbox{ for } i\leq j, \mbox{ then  }\alpha=0.
\label{induction}
\end{eqnarray}
This implies the proposition, because if $\alpha$ vanishes on every irreducible
component of the core then it vanishes on every irreducible projective subvariety $\overline{D}_i$ of the core.
The statement (\ref{induction})
then implies by
induction that $\alpha$ vanishes on the core.

To prove (\ref{induction})
consider the Mayer-Vietoris sequence of the covering
$ D_{\leq j}=D_{< j} \cup \overline{D}_j$.
\begin{eqnarray*} \dots \to H^k( D_{\leq j};\Q)\stackrel{\alpha}{\to} H^k(D_{< j};\Q)\oplus H^k(\overline{D}_j;\Q) \stackrel{\beta}{\to}
H^k(D_{<j}\cap
\overline{D}_j;\Q)\to \dots
\end{eqnarray*}
We show that the map  $\alpha$  is injective, which
will prove our claim. For this we show that
$ \beta$ is
surjective. This follows from the surjectivity of $H^k(\overline{D}_j; \Q) \to H^k(
\overline{D}_j\backslash D_j; \Q)$, because clearly
$D_{<j}\cap \overline{D}_j=\overline{D}_j\backslash D_j$.

To prove this we do Morse theory on the projective toric orbifold
$\overline{D}_j$. First it follows from Morse theory that
$H^*(\overline{D}_j;\Q)\to H^*(\overline{D}_j\backslash F_j;\Q)$ surjects. Moreover
we have that $ \overline{D}_j\backslash D_j$ is the core of the quasi-projective variety $\overline{D}_j\backslash F_j$. This means that $ \overline{D}_j\backslash D_j$ is
the set of points $x$ in $\overline{D}_j$ such that
$\lim_{z\to \infty} \lambda_v(z) x$ is not in $F_j$. Then the proof of
Theorem~\ref{circle} shows that $H^*(\overline{D}_j\backslash F_j;\Q)$ is
isomorphic with $H^*( \overline{D}_j\backslash D_j;\Q)$. This proves
(\ref{induction})
and in turn our Proposition~\ref{littlemma}.
\end{proof}

We finish this section with an explicit description of the cohomology ring of
$C(A,\theta)$, namely, we identify it with the cohomology of
the ambient semi-projective toric orbifold $X(A,\theta)$:

\begin{theorem}\label{circle}
The embedding of the core $C(A,\theta)$ in $X(A,\theta)$
induces  an isomorphism on cohomology with integer coefficients.
\end{theorem}

\begin{proof}
Let $v\in {\rm int}|\Sigma|$, $F_i$, $U_i$ and $D_i$ as above. We clearly have an
inclusion $D_{\leq j}\subset U_{\leq j}$.
We show by induction on $j$ that this inclusion induces an 
isomorphism on cohomology.
Consider the following commutative diagram:
$$\begin{array} {ccccccc}\dots \to & H^k(U_{\leq j},U_{<j};\Z)&
\to & H^k(U_{\leq j};\Z)&\to & H^k(U_{<j};\Z)&\to \dots\\
&\downarrow & & \downarrow && \downarrow & \\
\dots  \to & H^k(D_{\leq j},D_{<j};\Z)& \to & H^k(D_{\leq j};\Z)&\to & H^k(D_{<j};\Z)&\to \dots
\end{array}.
$$
The rows are the long exact sequence of the pairs $(U_{\leq j},U_{<j})$
and $(D_{\leq j},D_{<j})$ respectively.
The vertical arrows are induced by inclusion.
The last vertical arrow is an isomorphism by induction.

By excision
$H^k(U_{\leq j},U_{<j};\Z)\cong H^k(T(N_j), t_0 ;Z)$, where $N_j$ is the normal (orbi-)bundle to $U_j$ and
$T(N_j)$ is the Thom space $N_j\cup t_0$, where $t_0$ is the point at infinity.
Similarly $H^k(D_{\leq j},D_{<j};\Z)\cong H^k(T(D_j), t_0 ;\Z)$,
where $T(D_j)=D_{\leq j}/D_{<j}$ is the one point compactification of $D_j$, which is homeomorphic to the Thom space of $N_j|_{F_j}$,
 the negative bundle at $F_j$. Because $F_j$ is a
deformation retract of $U_j$ and because the normal bundle $N_j$ to $U_j$ 
in $U_{\leq j}$ restricts to
the normal bundle of $F_j$ in $D_j$, we find
that $T(D_j)$ is a deformation retract of $T(N_j)$.
Consequently the first vertical arrow is also an isomorphism. The
Five Lemma now delivers our assertion.
\end{proof}

\begin{remark} One can prove more, namely, that $C(A,\theta)$
is a deformation retract of $X(A,\theta)$.
This follows from Theorem \ref{circle}
and the analogous statement about the
fundamental group, which vanishes for both spaces.
Alternatively, one can use
Bott-Morse theory in the spirit
of the proof of \cite[Theorem 3.2]{Mi}
to get the homotopy equivalence.
\end{remark}

\section{Lawrence toric varieties}
\label{lawrence}

In this section we examine an important class of
toric varieties which are semi-projective but not projective.
We fix an integer $d \times n$-matrix $A$ as in  (\ref{ses}),
and we write  $A^\pm=[A,-A]$ for the $d \times 2n$-matrix
obtained by appending the negative of $A$ to $A$. The
corresponding vector configuration $\,\A^\pm=\A\cup -\A \,$
spans $\Z^d$ as a semigroup; in symbols, $\, \N\A^\pm = \Z\A = \Z^d$.
A vector $\theta$ is {\it generic}
with respect to $\A^{\pm}$ if it does not lie on
any hyperplane spanned by a subset of $\A$.

\begin{definition}
We call  $\,X(A^\pm,\theta)\,$ a {\em Lawrence toric variety},
for any generic vector $\theta \in \Z^d$.
\end{definition}

Our choice of name comes from the Lawrence construction
in polytope theory; see e.g.~Chapter 6 in \cite{Zi}.
The Gale dual of the centrally symmetric configuration
$\, \A^\pm \,$ is denoted $\Lambda(\B)$ and is called
the {\it Lawrence lifting} of $\B$. It consists
of $2n$ vectors which span $\Z^{2n-d}$.
The cone ${\rm pos} (\Lambda(\B))$ is the cone
over the $(2n-d-1)$-dimensional {\it Lawrence polytope}
with Gale transform $\A^\pm$.

Consider the even-dimensional affine space $\C^{2n}$ with
coordinates $z_1,\ldots,z_n,w_1,\ldots,w_n$.
We call a torus action on $\C^{2n}$ {\it symplectic}
if the products $\,z_1 w_1, \ldots, z_n w_n\,$
are fixed under this action.

\begin{proposition}
The following three classes of toric varieties coincide:
\begin{enumerate}
\item Lawrence toric varieties,
\item toric orbifolds which are GIT-quotients
of a symplectic torus action on  $\C^{2n}$ for some $n \in \N$,
\item toric varieties $X(\Sigma)$ where
$\Sigma$ is the cone over a regular triangulation of
a Lawrence polytope.
\end{enumerate}
\label{sympac}
\end{proposition}

\begin{proof}
This follows from Theorem~\ref{semi} using the observation
that a torus action on $\C^{2n}$ is symplectic if and only if it arises
from a matrix of the form  $A^\pm$. This means the action looks like
$$ z_i \, \mapsto \, t^{a_i} \cdot z_i \,\,, \qquad
 w_i \, \mapsto \, t^{- a_i} \cdot w_i \qquad \qquad (i = 1,2,\ldots,n) $$
Note that a polytope is Lawrence if and only
if its Gale transform is centrally symmetric.
\end{proof}

The matrix $A^\pm$ is unimodular if and only if the
smaller matrix $A$ is unimodular. Therefore
unimodularity of $A$ implies the  smoothness of the
Lawrence toric variety, by Corollary \ref{key}.
 An interesting feature of Lawrence toric varieties is
that the converse to this statement also holds:

\begin{proposition}
The Lawrence toric variety
$X(A^\pm,\theta)$ is smooth if and only if $A$ is unimodular.
\end{proposition}

\begin{proof}
The chamber complex $\Gamma(\A^\pm)$ is the arrangement
of hyperplanes spanned by subsets of $\A$. The vector $\theta$
is assumed to lie in an open cell of that arrangement.
For any column basis $C = \{a_{i_1},\ldots,a_{i_d}\}$
of the $d \times n$-matrix $\,A\,$ there exists a unique linear combination
$$
\lambda_1 a_{i_1} \, + \,
\lambda_2 a_{i_2} \, + \, \cdots \, + \,
\lambda_d a_{i_d} \quad = \quad \theta . $$
Here all the coefficients $\,\lambda_j \,$
are non-zero rational numbers.
We consider the polynomial ring
$$ \Z [ z,w ] \quad = \quad
\Z[x_1,\ldots,x_n, y_1,\ldots,y_n]. $$
The $2n$  variables are used to index the elements of
$\A^\pm$ and the elements of $\Lambda(\B)$. We set
$$\, \sigma(C,\theta) \quad = \quad
\{ \, x_{i_j} \, : \, \lambda_j > 0 \, \} \  \cup \,
\{ \, y_{i_j} \, : \, \lambda_j < 0 \, \}. $$
Its complement $\, \overline{\sigma}(C,\theta) \, = \,
\bigl\{x_1,\ldots,x_n, y_1,\ldots,y_n \bigr\} \backslash
\sigma(C,\theta)\,$
corresponds to a subset of $\Lambda(\B)$ which
forms a basis of $\R^{2n-d}$.
The triangulation $\Sigma_\theta$ of the Lawrence polytope defined by $\theta$
is identified with its set of maximal faces. This set equals
\begin{equation}
\label{lawtrig}
\Sigma_\theta \quad = \quad
\bigl\{ \, \overline{\sigma}(C,\theta)  \,\, \, : \,\,\,
C \,\, \hbox{is any column basis of}\,\, A \,\bigr\}.
\end{equation}
Hence the Lawrence toric variety $\, X(A^\pm,\theta) \, = \,
X(\Sigma_\theta)\,$ is smooth if and only if
every basis in $\Lambda(\B)$  spans the lattice $\Z^{2n-d}$
if and only if every column basis $C$ of $A$
spans $\Z^d$. The latter condition is equivalent to saying
that $A$ is a unimodular matrix.
\end{proof}

The $\Z^d$-graded polynomial ring  $\,{\Z}[x,y]$
is the homogeneous coordinate ring \cite{Cox} of
$X(\Sigma_\theta)$.

\begin{corollary}
The Stanley-Reisner ideal of the fan $\Sigma_\theta$ equals
\begin{equation}
\label{stanreis}
 I_\theta\quad = \quad \bigcap_{C} \,
\langle \sigma(C,\theta) \rangle
\quad \subset \,\,\, \Z[x,y],
\end{equation}
i.e. $I_\theta$ is
 the intersection of the monomial prime ideals generated by the sets
$\,\sigma(C,\theta) \,$ where $C$ runs over all column bases of $A$.
The irrelevant ideal of the
Lawrence toric variety $X(\Sigma_\theta)$  equals
\begin{equation}
\label{irrelevant}
 B_\theta \quad = \quad
\langle \, \prod  \sigma(C,\theta) \,: \, C \,\,
\hbox{is any column basis of} \,\, A \, \rangle
\quad \subset \,\,\, \Z[x,y].
\end{equation}
\end{corollary}

We now compute the cohomology of a Lawrence toric variety.
For simplicity of exposition we assume  $A$ is unimodular
so that $X(A^\pm,\theta)$ is smooth. The orbifold case
is analogous. First note
$$
{\rm Circ}(\A^\pm) \quad = \quad
\langle
 x_1 + y_1,
 x_2 + y_2,\ldots,
 x_n + x_n
\rangle \,\, + \,\,
{\rm Circ}(\A) , $$
where $\,{\rm Circ}(\A) \,$ is
generated by all linear forms
$\,\sum _{i=1}^n \lambda_i x_i  \,$ such that
$\,\lambda \, = \, (\lambda_1,\ldots,\lambda_n)\,$
lies in $\, {\rm ker}(A) = {\rm im}(B)$.
From Proposition~\ref{intcoho}, we have
$$ H^* \bigl( X(A^\pm, \theta);  \Z \,\bigr) \quad = \quad
\Z[x,y] \,/ \biggl(
\langle
 x_1 + y_1,
 x_2 + y_2,\ldots,
 x_n + y_n
\rangle \,\, + \,\,
{\rm Circ}(\A) \,\, + \,\,
I_\theta \biggr). $$
Let $\phi$ denote the $\Z$-algebra epimorphism
which collapses the variables pairwise:
$$ \phi \, : \,
\Z[x_1,\ldots,x_n,y_1,\ldots,y_n ] \, \rightarrow \, \Z[x_1,\ldots,x_n],
\quad x_i \mapsto x_i, \,\, y_i \mapsto - x_i \qquad
(i=1,2,\ldots,n). $$
Then we can rewrite the presentation of the cohomology ring as follows:
$$ H^* \bigl( X(A^\pm, \theta); \Z \bigr) \quad = \quad
\Z[x_1,\ldots,x_n] \,/ \biggl(
{\rm Circ}(\A) \,\, + \,\,
\phi(I_\theta) \biggr). $$
Clearly, the image of  the ideal  (\ref{stanreis})
under $\phi$ is the intersection of the ideals
$$ \phi \,\bigl( \, \langle \sigma(C,\theta) \rangle \, \bigr)
\quad = \quad
\langle \, x_i : \, i \in C \, \rangle
$$
where $C$ runs over the column bases of $A$.
Note that this ideal is independent of the choice of $\theta$
It depends only on $A$.
This ideal is called the {\it matroid ideal} of $\B$
and it is abbreviated by
\begin{eqnarray*}
M^*(\A) \quad &=& \quad \bigcap \bigl\{
\langle x_{i_1}, \ldots, x_{i_d} \rangle \,\,  : \,\,
\{a_{i_1} ,\ldots, a_{i_d} \} \subseteq \A
\quad \hbox{is linearly independent} \, \bigr\}
\\
&=& \langle x_{i_1}\cdots x_{i_k}: \{ b_{i_1},\dots, b_{i_k}\} \subseteq \B
\,\,\, \mbox{ is linearly dependent }\rangle \quad = \quad M(\B).
\end{eqnarray*}

We summarize what we have proved
concerning the cohomology of a Lawrence toric variety.

\begin{theorem}
\label{lawrencering}
The integral cohomology ring of a smooth Lawrence toric variety
$\, X(A^\pm, \theta) \,$ is independent of the choice of the generic vector
$\theta$ in $\Z^d$. It equals
\begin{equation}
 H^* \bigl( X(A^\pm, \theta); \Z \bigr) \quad = \quad
\Z[ x_1,\ldots,x_n] / \bigl(  {\rm Circ}(\A)  \, + \, M^*(\A) \,\bigr).
\end{equation}
The same holds for Lawrence toric orbifolds
with $\Z$ replaced by $\Q$.
\end{theorem}

\begin{remark} The independence of the cohomology ring on $\theta$ is
an unusual phenomenon in the GIT-construction. Usually,
the topology of the quotient changes when one crosses a wall.
Theorem \ref{lawrencering} says that this is
not the case for symplectic torus actions.
 An explanation of this fact
is offered through our Theorem~\ref{main},
as there are no walls in the hyperk\"ahler quotient construction.
\end{remark}

The ring $\,\Q[ x_1,\ldots,x_n]/ M^*(\A) \,$ is the
Stanley-Reisner ring of the matroid complex
(of linearly independent subsets) of
the $(n-d)$-dimensional configuration ${\cal B}$.
This ring is Cohen-Macaulay, and
$\, {\rm Circ}(\A)  \,$ provides a linear system
of parameters.  We write
$\,h(\B) = (h_0,h_1,\ldots,h_{n-d}) \,$ for its $h$-vector.
This is a well-studied quantity in combinatorics;
see e.g.~\cite{CC} and \cite[Section III.3]{St1}.

\begin{corollary}
\label{lawrencebetti}
The  Betti numbers of the Lawrence toric variety
$\, X(A^\pm, \theta)$
are independent of $\theta$, and they coincide with
the entries in the $h$-vector of the
rank $n-d$ matroid given by $\B$:
$$ {\rm dim}_\Q \, H^{2i}(X(A^\pm, \theta);\Q)
 \quad = \quad h_i(\B) .
\qquad \quad {\rm for} \quad i = 0,1,\ldots,n-d.$$
\end{corollary}

Our second result in this section concerns the core of a
Lawrence toric variety of dimension $2n-d$. We fix a generic vector $\theta$ in $\Z^d$.
The fan $\Sigma_\theta$ is the normal fan of the unbounded polyhedron
$$ P_\theta \quad = \quad
\bigl\{ \, (u,v) \in \R^{n} \oplus \R^n \,\, : \,\,
A u - A v \, = \, \theta , \,\, u,v \geq 0 \, \bigr\}. $$
As in the proof of Proposition~\ref{triangulation}, we chose any
vector $\psi \in \Z^n$ such that
$\, A \psi \, = \, - \theta $, and we consider the
following full-dimensional unbounded polyhedron in $\R^{2n-d}$:
$$ Q_\psi \quad = \quad
\bigl\{ \, (w,t) \in \R^{n-d} \oplus \R^n \, : \,
t \, \geq \, 0 \,, \,\,\,
B w \, + \, t \, \geq \, \psi \,\,\bigr\}. $$
The map $\,(w,t) \mapsto ( B w + t - \psi, t ) \,$ is an affine-linear
isomorphism from $Q_\psi$ onto $P_\theta$.

We define
${\mathcal H}(B,\psi)$ to be the arrangement of the following $n$
 hyperplanes in $\, \R^{n-d}$:
$$ \{ \, w \in \R^{n-d} \,\, : \,\, b_i \cdot w = \psi_i \, \}
\qquad \qquad (i=1,2,\ldots,n). $$
The arrangement ${\mathcal H}(B,\psi)$ is regarded
as a polyhedral subdivision of $\R^{n-d}$ into relatively
open polyhedra of various dimensions. The collection
of all such polyhedra which are bounded form a subcomplex,
called the {\it bounded complex} of ${\mathcal H}(B,\psi)$ and
denoted by ${\mathcal H}^{bd}(B,\psi)$.

\begin{theorem} The bounded complex ${\mathcal H}^{bd}(B,\psi)$
of the  hyperplane arrangement ${\mathcal H}(B,\psi)$
in $\R^{n-d}$ is isomorphic to the
complex of bounded faces of the $(2n-d)$-dimensional polyhedron
$\, Q_\psi \simeq P_\theta $.
\label{boundedbounded}
\end{theorem}

\begin{proof}
We define an injective map from $\R^{n-d}$ into
the polyhedron $\, Q_\psi \,$ as follows
\begin{equation}
\label{funnyembed}
   w \, \mapsto \,  \bigl(w, t \bigr),
\qquad \hbox{where} \,\,\,\,
t_i \, = \,
 {\rm max}\{0,\psi_i - b_i \cdot w \}.
\end{equation}
This map is linear on each cell of
the hyperplane arrangement ${\mathcal H}(B,\psi)$,
and the image of each cell is a face of $Q_\psi$.
In particular, every bounded cell of  ${\mathcal H}(B,\psi)$
is mapped to  a bounded face of $Q_\psi$ and
each unbounded cell of  ${\mathcal H}(B,\psi)$
is mapped to  an unbounded face of $Q_\psi$.
It remains to be shown that every bounded face of $Q_\psi$
lies in the image of the map  (\ref{funnyembed}).

Now, the image of   (\ref{funnyembed}) is the following
subcomplex in the boundary of our polyhedron:
\begin{eqnarray*}
& \{\, (w,t) \in Q_\psi \,\, : \,\,
t_i \cdot (b_i \cdot w + t_i - \psi_i ) \, = \, 0 \,\,\,
\hbox{for} \,\, i = 1,2,\ldots,n \, \bigr\}  \\
\simeq
& \{ \, (u,v) \in P_\theta \,\, : \,\,
u_i \cdot v_i
= \, 0 \,\,\,
\hbox{for} \,\, i = 1,2,\ldots,n \, \bigr\}
\end{eqnarray*}
Consider any face $F$ of $P_\theta$ which is not in this subcomplex,
and let $(u,v) $ be a point in the relative interior of $F$.
There exists  an index $i$ with
 $\, u_i > 0 \,$ and $\, v_i > 0 $.
Let $e_i$ denote the $i$-th unit
vector in $\R^n$. For every positive real $\lambda$,
the vector $\, (u + \lambda e_i, v + \lambda e_i) \,$
lies in $P_\theta$ and has the support as $(u,v)$.
Hence  $\, (u + \lambda e_i, v + \lambda e_i) \,$
lies in $F$ for all $\lambda \geq 0 $. This shows
that $F$  is unbounded.
\end{proof}

Theorem~\ref{boundedbounded} and Corollary \ref{semibetti} imply the following
enumerative result:

\begin{corollary} The Betti numbers of the Lawrence toric variety
$X(A^\pm,\theta)$ satisfy
$$ {\rm dim}_\Q \, H^{2i}(X(A^\pm, \theta);\Q) \quad = \quad
\sum_{i=k}^{n-d}(-1)^{i-k}
\binom{i}{k} f_i({\mathcal H}^{bd}(B,\psi)), $$ where $f_i({\mathcal H}^{bd}(B,\psi))$
denotes the number of $i$-dimensional
bounded regions in ${\mathcal H}(B,\psi)$.
\label{boundedbetti}
\end{corollary}

There are two natural geometric structures on any Lawrence  toric variety.
First the canonical bundle of $X(A^\pm,\theta)$
is trivial, because the vectors in $\A^\pm$ add to $0$. This
means that $X(A^\pm,\theta)$ is a {\em Calabi-Yau variety}.
Moreover, since the symplectic $\T^d_\C$-action
preserves the natural Poisson structure on $\C^{2n}\cong \C^n\oplus (\C^n)^*$,
the GIT quotient $X(A^\pm,\theta)$ inherits a natural {\em holomorphic Poisson structure}.
The holomorphic symplectic leaves of this Poisson structure are what we call
{\em toric hyperk\"ahler manifolds}.
The special leaf which   contains the
core of $X(A^\pm, \theta)$ will be called the
{\em toric hyperk\"ahler variety}. We present these definitions
in complete detail in the following two sections.

\section{Hyperk\"ahler quotients}
\label{hkquotients}

Our aim is to  describe an algebraic approach to the toric
hyperk\"ahler manifolds of Bielawski and Dancer \cite{BD}.
In this section we sketch the original
differential geometric construction in \cite{BD}. This construction
is the hyperk\"ahler
analogue to the construction of toric varieties using K\"ahler quotients.
We first briefly review the latter.
Fix the   standard Euclidean bilinear form on  $\C^n$,
$$g(z,w) \quad = \quad\sum_{i=1}^n \left(
{\rm re}( z_i) {\rm re}( w_i) + {\rm im}( z_i) {\rm im}( w_i)\right).$$
The corresponding
{\em K\"ahler form}  is
$$ \omega(z,w) \quad = \quad g(iz,w) \quad = \quad \sum_{i=1}^n \left(
{\rm re}( z_i) {\rm im}( w_i)
- {\rm im}( z_i) {\rm re} (w_i)\right).$$
Let $A$ be as in (\ref{ses}) and consider the real torus $\T^d_\R$
which is the maximal compact subgroup of $\T^d_\C$. The group
$\T^d_\R$  acts on $\C^n$
preserving the K\"ahler structure.
This action has the moment map
\begin{equation}
\label{kmomentmap}
\mu_\R \,:\,  \C^n  \,\to\,  ({\mathfrak t}_\R^d)^*\cong \R^d \, , \quad
 (z_1,\dots,z_n)
 \,\mapsto \, {1\over 2} \sum^n_{i=1} |z_i|^2 a_i.
\end{equation}
Fix $\xi_\R\in \R^d$. The {\em K\"ahler quotient}
 $\,X(A,\xi_\R)=\C^n//_{\xi_\R} \T^d_\R=\mu_\R^{-1}(\xi_\R)/\T^d_\R\,$
inherits a K\"ahler structure from $\C^n$ at its smooth points.
If $\xi_\R=\theta$ lies in the lattice $ \Z^d $
then there is a  biholomorphism
between the smooth loci in the GIT quotient $X(A,\theta)$ and the K\"ahler
quotient $X(A,\xi_\R)$. Hence if $A$ is unimodular and $\theta$ generic
then the complex manifolds $X(A,\theta)$ and $X(A,\xi_\R)$ are biholomorphic.

Now we turn to toric hyperk\"ahler manifolds.
Let $\H$ be the skew field of {\em quaternions}, the
$4$-dimensional real vector space
with  basis $1,i,j,k$  and associative algebra structure
given by  $\,i^{2}=j^{2}=k^{2}=ijk=-1$. Left
multiplication by $i$ (resp. $j$ and $k$)
defines complex structures $I:\H\to \H$, with $I^2=-{\rm Id}_\H$,
(resp. $J$ and $K$) on $\H$. We now put the
flat metric $g$ on $\H$ arising from
the standard Euclidean scalar product on
$\H\cong \R^4$
with $1,i,j,k$ as an orthonormal basis.
This is called a {\em hyperk\"ahler metric}  because it is
a K\"ahler metric with respect to
all three complex structures $I$, $J$ and $K$.
It means that the differential $2$-forms, the so-called {\em K\"ahler forms},
given by $\omega_I(X,Y)=g(IX,Y)$ for tangent vectors $X$ and $Y$,
and the analogously defined $\omega_J$ and
$\omega_K$ are closed.
A special orthogonal transformation,
with respect to this metric, is
said to preserve the hyperk\"ahler structure if
it commutes with all three complex structures $I$, $J$ and $K$ or equivalently if it preserves
the K\"ahler forms $\omega_I$, $\omega_J$ and $\omega_K$.
The group of such transformations,
the unitary symplectic group $Sp(1)$, is
generated by multiplication by unit quaternions from the right.
A maximal abelian subgroup $\T_\R\cong U(1)\subset Sp(1)$
is thus specified by a choice of a unit quaternion.
We break the symmetry between $I$, $J$ and $K$ and
choose the maximal
torus generated by multiplication from the right by the unit quaternion $i$.
Thus $U(1)$ acts on
$\H$ by sending $\xi$ to $\xi \exp(\phi i) $,
 for $exp(\phi i) \in U(1)\subset \R\oplus\R i \cong \C$.
It follows from (\ref{kmomentmap})  that
the moment map $\mu_I: \H\to \R$ with respect to the symplectic form
$\omega_I$ is given by
\begin{equation}
\label{plusplusminusminus}
 \mu_I(x+yi+uj+vk) \quad = \quad
\mu_I(x+yi+(-ui+v)k)
\quad = \quad {1\over 2} (x^2+y^2-u^2-v^2).
\end{equation}
Similarly we obtain
formulas for $\mu_J$ and $\mu_K$ by writing
down the eigenspace decomposition in the respective complex
structures:
\begin{eqnarray*}
\mu_J(x+yi+uj+vk) &=& \mu_J\left[\left({y+u\over \sqrt{2}}+{-x-v\over \sqrt{2}} j\right){i+j\over \sqrt{2}}+
\left({y-u\over \sqrt{2}}j + {-x+v\over \sqrt{2}} \right){k-1\over \sqrt{2}}\right]\\&=& yu+xv, \\
\mu_K(x+yi+uj+vk) &=& \mu_K\left[\left({y+v\over \sqrt{2}}+{-x+u\over \sqrt{2}} k\right){i+k\over \sqrt{2}}+
\left({y-v\over \sqrt{2}}+{x+u\over \sqrt{2}} k\right){i-k\over \sqrt{2}}\right]\\ &=& yv-xu.\end{eqnarray*}
We now consider the map $\, \mu_\C=\mu_J+i\mu_K \,$ from $\H$ to $\C$.
It can be thought of
as the holomorphic moment map for the $I$-holomorphic action of
$\T_\C\supset \T_\R$ on $\H$ with respect to the $I$-holomorphic symplectic
form $\omega_\C=\omega_J+i \omega_K$.
If we identify  $\H$ with $\C \oplus \C$ by introducing
two complex coordinates, $z=x+iy\in \R\oplus \R i\cong \C$ and
$w= v-ui \in \R\oplus \R i\cong \C$, then
the $I$-holomorphic moment map $\mu_\C : \H \rightarrow \C$ is given
algebraically by multiplying complex numbers:
\begin{equation}
\label{remarksimp}
\mu_\C(z,w) \quad = \quad \mu_J(z,w) \, + \, i\mu_K(z,w)
\quad = \quad yu+xv \, + \, i(yv-xu)
\quad = \quad zw .
\end{equation}

The discussion in the previous paragraph
generalizes in an obvious manner to
$\H^n$ for $n > 1$. Indeed, the
$n$-dimensional quaternionic space $\H^n$
has three complex structures
$I$,$J$ and $K$, given by left multiplication with $i,j,k\in \H$.
Putting the flat metric
$g_n=g^{\oplus n}$ on $\H^n$ yields a hyperk\"ahler metric,
i.e.~the differential $2$-forms
$\omega_I(X,Y)=g_n(IX,Y)$ and similarly $\omega_J$ and $\omega_K$
are K\"ahler (meaning closed) forms.
The automorphism group of this hyperk\"ahler structure is
the unitary symplectic group $Sp(n)$.
We fix the maximal torus $\T_\R^n=U(1)^n\subset
Sp(n)$ given by the following definition.
For $\lambda=(\exp(\phi_1 i ),
\exp(\phi_2 i),\dots, \exp(\phi_n i))\in \T_\R^n$ and
$(\xi_1,\xi_2,\dots,\xi_n)\in \H^n$ we set
\begin{equation}
\label{nnnaction}
\lambda (\xi_1,\xi_2,\dots,\xi_n)
\quad = \quad  (\xi_1\exp(\phi_1 i ),
 \xi_2 \exp(\phi_2 i),\dots,\xi_n \exp(\phi_n i)).
\end{equation}
As in the $n=1$ case above, this fixes an isomorphism
$\, \H^n\cong \C^n\oplus \C^n \,$
where two complex vectors $z,w\in \C^n\cong \R^n\oplus i\R^n $
represent the quaternionic vector
$z+wk\in \H^n \cong \R^n\oplus i\R^n  \oplus j\R^n  \oplus k\R^n$.
Expressing vectors in $\H^n$ in these complex coordinates,
the torus  action  (\ref{nnnaction}) translates into
\begin{eqnarray}
\label{action}\lambda  (z,w) \quad = \quad(\lambda z, \lambda^{-1} w)
\qquad \hbox{for} \,\,\, \lambda\in \T_\R^n
\quad \hbox{and} \quad (z,w) \in \H^n.
\end{eqnarray}

The  toric hyperk\"ahler manifolds in  \cite{BD}
are constructed by choosing a subtorus
$\T_\R^d\subset \T_\R^n$ and taking the hyperk\"ahler quotient
\cite{HKLR} of $\H^n$ by $\T_\R^d$.
We do this by choosing integer matrices  $A$ and $B$ as in (\ref{ses})
and  (\ref{contraC}).
The subtorus $\T^d_\R $ of $ \T^n_\R$
acts on $\H^n$ by (\ref{action}) preserving
the hyperk\"ahler structure.
The {\em hyperk\"ahler moment map} of the action (\ref{action}) of
$\T_\R^d$ on $\H^n$ is defined by
$$\mu=(\mu_I,\mu_J,\mu_K) \, : \, \H^n \, \to  \,
({\mathfrak t}_\R^d)^*\otimes \R^3,$$ where $\mu_I$, $\mu_J$ and
$\mu_K$ are the K\"ahler moment maps with respect to $\omega_I$, $\omega_J$ and $\omega_K$ respectively.
Using the formulas (\ref{plusplusminusminus}) and
(\ref{remarksimp}), the components
of $\mu$ are in complex coordinates as follows:
\begin{equation}
\mu_\R(z,w) \,\,\, := \,\,\, \mu_I(z,w)
\quad = \quad {1 \over 2} \sum_{i=1}^n
(|z_i|^2-|w_i|^2) \cdot a_i
\quad \in \,\,\, ({\mathfrak t}_\R^d)^*, \qquad
\label{realmm}\end{equation}
\begin{equation}\mu_\C(z,w) \,\,\, :=  \,\,\, \mu_J(z,w)+i
\mu_K(z,w) \quad = \quad
\sum_{i=1}^n z_i w_i \cdot a_i
\quad \in \,\,\, ({\mathfrak t}^d_\R)^*\otimes
\C\cong({\mathfrak t}^d_\C)^*. \label{holmm1}\end{equation}
Here $a_i$ is the $i$-th column vector of the matrix $A$.
We can also think of
$\mu_\C$ as the  moment map for the $I$-holomorphic action
of $\T^d_\C$ on $\H^n$ with respect to  $\,\omega_\C=\omega_I+i\omega_K$.
Now take
$$ \xi=(\xi^1,\xi^2,\xi^3)
\,\,\, \in \,\,\, ({\mathfrak t}_\R^d)^* \otimes \R^3$$ and introduce
$\xi_\R=\xi^1\in ({\mathfrak t}_\R^d)^*$ and $\xi_\C=\xi^2+i\xi^3\in ({\mathfrak t}_\C^d)^*$ so we can write
$\xi=(\xi_\R,\xi_\C)\in  ({\mathfrak t}_\R^d)^* \oplus ({\mathfrak t}_\C^d)^*$.
The {\em hyperk\"ahler quotient}
of $\H^n$ by the action (\ref{action}) of the torus
$\T_\R^d$ at level $\xi$ is defined as
\begin{eqnarray}\label{quotient} Y(A,\xi) \,\,\, := \,\,\,
\H^n /\!/ \! / \! /_{\xi} \T^d_\R \quad
:= \quad \mu^{-1}(\xi)/\T_\R^d \quad = \quad
\left( \mu_\R^{-1}(\xi_\R) \cap
\mu_\C^{-1}(\xi_\C) \right)/ \T_\R^d.\end{eqnarray}
By a theorem of \cite{HKLR}, this quotient has a canonical
hyperk\"ahler structure on its smooth locus.

Bielawski and Dancer show in \cite{BD}
that if $\xi\in ({\mathfrak t}_\R^d)^*\otimes\R^3 $
is generic then $Y(A,\xi)$ is an orbifold,
and it is smooth if and only if $A$ is unimodular.
Since $\xi$ is generic outside a set of codimension
three in  $({\mathfrak t}_\R^{d})^*\otimes\R^3 $, they can show
that the topology and therefore the cohomology of the toric hyperk\"ahler
manifold is independent on $\xi$. In what follows we consider
vectors $\xi$ for which
$\xi_\C=0 $ in $\C^d$  and
$\xi_\R=\theta \in \Z^d\subset \R^d\cong ({\mathfrak t}_\R^d)^* $.
The underlying
complex manifold in complex structure $I$ of the
hyperk\"ahler manifold $Y(A,(\theta,0_\C))$
has a purely algebraic description as explained in
the next section.

\section{Algebraic construction of toric hyperk\"ahler varieties}
\label{konno}

The $\Z^d$-graded polynomial
ring $\C[z,w] = \C[z_1,\dots,z_n,w_1,\dots,w_n]$, with
the grading given by $A^{\pm}=[A,-A]$, is
the homogeneous coordinate ring of
the  Lawrence toric variety $X(A^\pm,\theta)$.
By a result of Cox \cite{Cox},
closed subschemes of $X(A^\pm,\theta)$
correspond to homogeneous ideals in $\C[z,w]$
which are saturated with respect to the
irrelevant ideal $B_\theta$ in (\ref{irrelevant}).
Let us now consider the ideal
\begin{eqnarray}
\label{cocircuitideal}
 \Cocirc \,\,\, := \,\,\,
\langle \,\sum_{i=1}^n a_{ij} z_i w_i  \,\,| \,\, j=1,\dots,d \,
\rangle \quad \subset \quad \C[z,w],
\end{eqnarray}
whose generators are  the components of the holomorphic moment map
$\mu_\C$ of (\ref{holmm1}). The ideal $\Cocirc$
is clearly homogeneous and it is a complete intersection.
We assume that none of the row vectors of the matrix $B$ is zero.
Under this hypothesis, the ideal  $\Cocirc$ is a prime ideal.

\begin{definition}
\label{complete} $\!\!\!\!$
 The {\em toric hyperk\"ahler variety}
$Y(A,\theta)$ is the irreducible subvariety of
the Lawrence toric variety  $X(A,\theta)$
defined by the homogeneous ideal $\Cocirc$
in the coordinate ring $\C[z,w]$ of $X(A,\theta)$.
\end{definition}

\begin{proposition} If $\theta$ is generic then
the toric hyperk\"ahler variety $Y(A,\theta)$ is an orbifold.
It is smooth if and only if the matrix $A$ is unimodular.
\end{proposition}

\begin{proof} It follows from (\ref{holmm1})
 that a point in $\C^{2n}$ has a finite stabilizer
under the group $\T_\C^d$ if and only if the point is regular for $\mu_\C$ of
(\ref{holmm1}), i.e. if the derivative of $\mu_\C$ is surjective there.
This implies that, for $\theta$ generic, the toric hyperk\"ahler variety
$Y(A,\theta)$ is an orbifold because then the variety
$X(A^\pm,\theta)$ is an orbifold.
For the second statement note that
if $A$ is unimodular then   $X(A^\pm,\theta)$ is smooth, consequently
$Y(A,\theta)$ is also smooth. However, if $A$ is not unimodular
then $X(A^\pm,\theta)$ has
orbifold singularities which lie in the core.
Now the core $C(A^\pm,\theta)$
lies entirely  in $Y(A,\theta)$,
by Lemma~\ref{cores} below, thus $Y(A,\theta)$ inherits singular
points from $X(A^\pm,\theta)$.
\end{proof}

We can now prove that our toric hyperk\"ahler varieties are biholomorphic to
the toric hyperk\"ahler manifolds of the previous section.

\begin{theorem} Let $\xi_\R=\theta\in \Z^d\subset (\t_\R^d)^*\cong \R^d$
for generic $\theta$. Then the
toric hyperk\"ahler manifold $Y(A,(\xi_\R,0))$ with complex
structure $I$ is biholomorphic
with the toric hyperk\"ahler variety $Y(A,\theta)$.
\end{theorem}

\begin{proof} Suppose $A$ is unimodular. The general theory of
K\"ahler quotients (e.g.~in \cite{Ki})
implies that the Lawrence toric variety
$X(A^\pm,\theta)$ and the corresponding K\"ahler quotient
$X(A^\pm,\xi_\R)=\mu_\R^{-1}(\xi_\R)/ \T_\R^d$ are  biholomorphic,
where $\mu_\R$ is defined in (\ref{realmm})
and $\xi_\R=\theta\in \Z^d\subset \R^d\cong (\t_\R^d)^*$.
Now the point is
that $\mu_\C:\H^n\to \C^d$ is invariant under the action of
$\T_\R^d$ and therefore descends to a map
on $X(A^\pm,\xi_\R)=\mu_\R^{-1}(\xi_\R)/ \T_\R^d$ and similarly on $X(A^\pm,\theta)$ making the following
diagram commutative: \begin{eqnarray*} \begin{array}{cccc}\mu^\xi_\C :&X(A^\pm,\xi_\R) &\to &\C^d \\ &\cong &&\cong\\
 \mu^\theta_\C :& X(A^\pm,\theta) & \to &\C^d \end{array}.
\end{eqnarray*} It follows that
$Y(A,(\xi_\R,0))=(\mu_\C^\xi)^{-1}(0)$ and $Y(A,\theta)=(\mu_\C^\theta)^{-1} (0)$
are biholomorphic as claimed.

The proof is similar in the case when the spaces have orbifold
singularities.
\end{proof}

Recall the affinization map
$\pi_X: X(A^\pm,\theta)\to X(A^\pm,0)$ from (\ref{pix}),
and the analogous map $\pi_Y: Y(A,\theta)\to Y(A,0)$. These fit together
in the following commutative diagram:
 \begin{eqnarray*}\begin{array}{ccc}  Y(A,\theta)&\stackrel{\pi_Y}{\to}& Y(A,0)\\
i_{\theta}\downarrow&&\downarrow i_{0}\\ X(A^\pm,\theta)&\stackrel{\pi_X}{\to}& X(A^\pm,0)\\
\mu_\C^\theta\downarrow &&\downarrow \mu_\C^0 \\ \C^d&\cong & \C^d\end{array}, \end{eqnarray*}
where $i_\theta:Y(A,\theta)\to X(A^\pm,\theta)$
denotes the natural embedding in Definition~\ref{complete}
by the preimage of $\mu_\C^\theta$ at
$0\in \C^d$. From this we deduce the following lemma:

\begin{lemma}
The cores of the Lawrence
toric variety and of the toric hyperk\"ahler variety coincide,
that is,
$\, C(A^\pm,\theta)=\pi_X^{-1}(0)= \pi_Y^{-1}(0)$.
\label{cores} \end{lemma}

\begin{remark} It is shown in \cite{BD} that the core of the toric hyperk\"ahler manifold
$Y(A,\theta)$ is the preimage of the bounded complex in the hyperplane arrangement ${\mathcal H}(\B,\psi)$
by the hyperk\"ahler moment map. We know from Theorem~\ref{coredescr} that the core of the Lawrence toric
variety equals the preimage of $P^{bd}_\theta$ under the K\"ahler moment map. 
Thus Theorem~\ref{boundedbounded} is a combinatorial analogue of
Lemma~\ref{cores}.
\end{remark}

We need one last ingredient in order to prove
the theorem stated in the Introduction.

\begin{lemma}
\label{isocohommm}
The embedding of the core $C(A,\theta)$ in $Y(A,\theta)$
gives an isomorphism in cohomology.
\end{lemma}

\begin{proof}
Consider the $\T_\C$-action on the Lawrence toric variety $X(A^\pm,\theta)$
defined by the vector $v=\sum_{i=1}^n b_i\in \Z^{n-d}$. This action
comes from multiplication by non-zero complex numbers on the vector space
$\C^{2n}$.  The holomorphic moment map
$\mu_\C$ of (\ref{holmm1}) is homogeneous with respect to multiplication
by a non-zero complex number, and consequently  $\mu_\C^\theta$
 is also homogeneous with respect to the circle action
$\lambda_v$. It follows that this $\T_\C$-action leaves the toric
hyperk\"ahler variety invariant. Moreover, since  $v$ is in
the interior of ${\rm pos}({\B})$, all the results in
Section \ref{corevariety} are valid for this $\T_\C$-action on
$X(A^\pm,\theta)$. Now the proof of Theorem~\ref{circle} can be repeated
verbatim to
show that the cohomology of $Y(A,\theta)$
agrees with the cohomology of the core.
\end{proof}

\noindent {\sl Proof of Theorem~\ref{main}}:
1.= 3. is a consequence of Lemma~\ref{cores} and
Lemma~\ref{isocohommm}.

2.= 3. This is a consequence of Theorem~\ref{circle}.

1.= 4. is the content of Theorem~\ref{lawrencering}.
 $\square$

\begin{remark} 1. In fact, we could claim more
than the isomorphism of cohomology rings in
Theorem~\ref{main}.
The remark after
Theorem~\ref{circle} implies
that the spaces $\, C(A^\pm,\theta)\subset
Y(A,\theta) \subset X(A^\pm, \theta)\,$ are
deformation retracts in one another. A similar result appears in
 \cite[Theorem 6.5]{BD}.

2. The result 2.=4. in the smooth case was proven by Konno in \cite{Ko}.

3. We deduce from Theorem~\ref{main}, Corollary~\ref{lawrencebetti}
and Corollary~\ref{boundedbetti} the following 
formulas for Betti numbers.
The second formula is due to 
Bielawski and Dancer  \cite[Theorem 6.7]{BD}.
\end{remark}

\begin{corollary}
\label{hyperbetti}
The Betti numbers of the toric hyperk\"ahler variety
$Y(A,\theta)$ agree with:
\begin{itemize}
\item the $h$-numbers of the matroid of $\B$: $
\quad b_{2k}(Y(A,\theta))\,= \, h_k(\B)$.
\item the following linear combination of the number of bounded regions of the affine hyperplane arrangement ${\mathcal H}(B,\psi)$: 
\begin{equation}
\label{alternating}
 b_{2k}(Y(A,\theta)) \quad = \quad \sum_{i=k}^{n-d}(-1)^{i-k}
\binom{i}{k} f_i({\mathcal H}^{bd}(B,\psi)). 
\end{equation}
\end{itemize}
\end{corollary}

This corollary shows the importance
of the combinatorics of the bounded complex ${\mathcal H}^{bd}(B,\psi)$ in the topology
of $Y(A,\theta)$ and $X(A^\pm,\theta)$. This intriguing connection will be more apparent 
in the next section. Before we get there we infer some important properties of 
the bounded complex from Corollary 9.1 and Theorem E of \cite{Za}.

\begin{proposition}\label{isispure} The bounded complex 
${\mathcal H}^{bd}(B,\psi)$ is pure-dimensional.
If $\psi$ is generic and $\B$ is coloop-free then
every maximal face of ${\mathcal H}^{bd}(B,\psi)$
is an $(n-d)$-dimensional simple polytope.
\end{proposition}

A {\em coloop} of $\B$ is a vector $b_i$
which lies in every column basis of $B$.
This is equivalent to $a_i$ being zero.
Note that if $A$ has a zero column then
we can delete it to get $A^\prime$, which  
 means that $Y(A,\theta)=Y(A^\prime,\theta)\times \C^{2}$ and similarly for the Lawrence toric 
variety. Therefore we will assume in the next section 
that none of the columns of $A$ is zero.

\section{Cogenerators of the cohomology ring}
\label{cogenerators}

There are three natural presentations of the
cohomology ring of the toric hyperk\"ahler variety
$Y(A,\theta)$ associated with a $d \times n$-matrix $A$
and a generic vector $\theta \in \Z^d$.
In these presentations $\,H^*(Y(A,\theta);\Q)\,$ is
expressed as a  quotient of
the polynomial ring  $\,  \Q[x,y] \,$ in $2n$ variables,
as a quotient of
the polynomial ring $\,  \Q[x] \,$ in $n$ variables,
or as a quotient
 of the polynomial ring $\,\Q[\,t \,] \,\simeq \,
 \Q[x]/\Circ \,$ in $d$ variables,
respectively. In this section we compute systems
of cogenerators for  $\,H^*(Y(A,\theta);\Q)\,$ relative
to each of the three presentations. As an application we show that
the  Hard Lefschetz Theorem holds for
toric hyperk\"ahler varieties, and we
discuss some implications for the combinatorial problem
of classifying  the $h$-vectors of matroid complexes.

We begin by reviewing the definition of cogenerators
of a homogeneous polynomial ideal. Consider the commutative
polynomial ring generated by a basis of derivations on affine $m$-space:
$$ \Q[\partial] \quad = \quad \Q[ \partial_1,\partial_2,\ldots,\partial_m]. $$
The polynomials in $\Q[\partial]$ act as linear
differential operators with constant coefficients on
$$ \Q[x] \quad = \quad \Q[ x_1, x_2, \ldots, x_m ]. $$
If $\Gamma$ is any subset of $\Q[x]$ then its
{\it annihilator}  $\, {\rm Ann}(\Gamma) \,$ is
the ideal in $\Q[\partial]$ consisting of all
linear differential operators with constant coefficients
which annihilate all polynomials in $\Gamma$.
If $I$ is any zero-dimensional homogeneous ideal in
$\Q[\partial]$ then there exists a finite set
$\Gamma$ of homogeneous polynomials in $\Q[x]$
such that $\, I \, = \, {\rm Ann}(\Gamma)$.
We say that $\Gamma$ is a set of {\it cogenerators} of $I$.
If $\Gamma$ is a singleton, say, $\Gamma = \{p\}$,
then  $\, I \, = \, {\rm Ann}(\Gamma)\,$ is a {\it Gorenstein ideal}.
In this case, the  polynomial $p = p(x)$ which cogenerates
$I$ is unique up to scaling. More generally, if
all polynomials in $\Gamma$ are homogeneous
of the same degree then  $\, I \, = \, {\rm Ann}(\Gamma)\,$ is a
{\it level ideal}. In this case, the
$\Q$-vector space spanned by $\Gamma$ is unique,
and it is desirable for $\Gamma$ to be
a nice basis for this space.

We replace the vector $\,\psi = (\psi_1,\ldots,\psi_n)\,$ in
Theorem~\ref{boundedbounded} by an indeterminate vector $\,x  = (x_1,\ldots,x_n)$
which ranges over a small neighborhood of $\psi$ in $\R^n$.
For $x$ in this neighborhood,
the polyhedron  $Q_x$ remains simple and
combinatorially isomorphic to  $ Q_\psi $,
and the hyperplane arrangement ${\mathcal H}(B,x)$
remains  isomorphic to ${\mathcal H}(B,\psi)$.
Let $\,\Delta_1,\ldots,\Delta_r \,$ denote the
maximal bounded regions of  ${\mathcal H}(B,x)$. These are
$(n-d)$-dimensional simple polytopes, by Proposition \ref{isispure}
and our assumption that $\B$ is coloop-free.
They can be identified with the maximal bounded faces
of the $(2n-d)$-dimensional polyhedron $Q_x$,
by Theorem~\ref{boundedbounded}.
The volume of the polytope $\Delta_i$
is a  homogeneous polynomial in $x$ of degree $n-d$
 denoted
 $$ V_i(x) \,\,\, = \,\,\, V_i(x_1,\ldots,x_n)
\quad = \quad {\rm vol}(\Delta_i)
\qquad (i=1,2,\ldots,r) $$

\begin{theorem}
\label{hypercogenerator}
 The volume polynomials $V_1,\ldots,V_r$
form a basis of cogenerators for the cohomology ring
of the Lawrence toric variety  $X(A^\pm,\theta)$ and
of the toric hyperk\"ahler variety  $Y(A,\theta)$:
\begin{equation}
\label{hypercogennn}
\,H^*(Y(A,\theta);\Q) \quad = \quad
\Q[ \partial_1,\partial_2,\ldots,\partial_n]/{\rm Ann}
\bigl( \{V_1,V_2,\ldots,V_r\}\bigr).
\end{equation}
\end{theorem}

\begin{proof} Each simple polytope $\Delta_i$ represents an $(n-d)$-dimensional
projective toric variety $X_i$. The core
$C(A^\pm,\theta)$ is glued from the toric varieties $X_1,\ldots,X_r$,
and it has the same cohomology as $X(A^\pm,\theta)$ and
$Y(A,\theta)$ as proved in Theorem~\ref{main}.
Hence we get a natural ring epimorphism induced
from the inclusion of each
toric variety $X_i$ into  the core $C(A^\pm,\theta)$:
\begin{equation}
\label{topphi}
\phi_i \quad : \quad H^* \bigl( C(A^\pm,\theta) ;\Q \bigr)
\,\,\rightarrow \,\, H^*(X_i; \Q).
\end{equation}
In terms of coordinates, the map $\phi_i$ is described as follows:
\begin{equation}
\label{algphi}
\phi_i \,\,\, : \,\,\,
\Q[\partial_1,,\ldots,\partial_n]/ \bigl( M(\B) + {\rm Circ}(\A) \bigr)
\,\, \rightarrow\,\,
\Q[\partial_1,,\ldots,\partial_n]/ \bigl( I_{\Delta_i} + {\rm Circ}(\A) \bigr),
\end{equation}
where $I_{\Delta_i}$ is the Stanley-Reisner ring of the simplicial
normal fan of the polytope $\Delta_i$. Each facet of $\Delta_i$ has
the form $\, \{\,w \in \Delta_i \, : \, b_j \cdot w \, = \, \psi_j \,\}\,$
for some $j \in \{1,2,\ldots,n \}$. The ideal $I_{\Delta_i}$ is generated by
all monomials $\,\partial_{j_1} \partial_{j_2} \cdots \partial_{j_s}\,$
such that the intersection of the facets
$\, \{\,w \in \Delta_i \, : \, b_{j_\nu} \cdot w \, = \, \psi_{j_\nu} \,\}$,
for $\nu = 1,2,\ldots,s$, is the empty set.
By the genericity hypothesis on $\psi$, this will happen
if $\,\{ b_{j_1}, b_{j_2},\ldots, b_{j_s}\} \,$ is linearly dependent,
or, equivalently, if  $\,\partial_{j_1} \partial_{j_2} \cdots \partial_{j_s}\,$
lies in the  matroid ideal $M(\B)$. We conclude that
$\, M (\B) \subseteq I_{\Delta_i}$, and the map $\phi_i$ in
(\ref{algphi}) is induced by this inclusion.

Proposition~\ref{littlemma} implies that
\begin{equation}
\label{crucialclaim}
{\rm ker}(\phi_1) \, \cap \,
{\rm ker}(\phi_2) \, \cap \,\, \ldots \,\,
\, \cap \, {\rm ker}(\phi_r) \quad = \quad \{0\}.
\end{equation}
Here is an alternative
proof for this in the toric hyperk\"ahler case.
We first note that the top-dimensional cohomology of
an equidimensional union of projective varieties equals the direct sum of the pieces:
\begin{equation}
\label{pieces}
H^{2n-2d} \bigl( C(A^\pm,\theta) ;\Q \bigr) \quad \simeq \quad
H^{2n-2d} \bigl( X_1 ;\Q \bigr) \, \oplus \,
\cdots \,\oplus \,
H^{2n-2d} \bigl( X_r ;\Q \bigr) ,
\end{equation}
and the restriction of the map $\phi_i$ to
degree $2n-2d$ is the $i$-th coordinate projection
in this direct sum. In particular, (\ref{crucialclaim})
holds in the top degree. We now use a theorem of
Stanley \cite[Theorem III.3.4]{St1} which states
that the Stanley-Reisner ring of a matroid is level.
Using condition (j) in \cite[Proposition III.3.2]{St1},
this implies that the socle of our cohomology ring
$\, H^{*} \bigl( C(A^\pm,\theta) ;\Q \bigr) \,$
consists precisely of the elements of degree $2n-2d$.
Suppose that (\ref{crucialclaim}) does not hold,
and pick a non-zero element $p(\partial)$ of maximal degree in the
left hand side. The cohomological degree of $p(\partial)$ is strictly
less than $2n-2d$ by (\ref{pieces}). For any generator $\partial_j$ of
$H^{*} \bigl( C(A^\pm,\theta) ;\Q \bigr) $,  the product
$\partial_j \cdot p(\partial)$ lies in the left hand side of  (\ref{crucialclaim})
because $\,\phi_i(  \partial_j \cdot p(\partial)) = \phi_i(  \partial_j)
\cdot \phi( p(\partial)) = 0 $. By the maximality hypothesis
in the choice of $p(\partial)$, we conclude that
$\partial_j \cdot p(\partial) = 0$ in $H^{*} \bigl( C(A^\pm,\theta) ;\Q \bigr) $
for all $j = 1,2,\ldots,n$. Hence $p(\partial)$ lies in the socle
of  $H^{*} \bigl( C(A^\pm,\theta) ;\Q \bigr) $.
By Stanley's Theorem, this means that $p(\partial)$ has
cohomological degree $\,2n-2d $. This is a contradiction and our claim follows.

The result (\ref{crucialclaim}) which we just proved  translates into the
following ideal-theoretic statement:
\begin{equation}
\label{interseccc}
 M(\B) \, + \, {\rm Circ}(\A) \quad = \quad
\bigcap_{i=1}^r \, \bigl( \, I_{\Delta_i} \, + \, {\rm Circ}(\A) \,\bigr).
\end{equation}

Since $X_i$ is a projective orbifold, the ring $\,H^*(X_i;\Q)\,$ is a Gorenstein ring.
A result of Khovanskii and Pukhlikov \cite{KP} states that its cogenerator
is the volume polynomial, i.e.
$$ I_{\Delta_i} \, + \, {\rm Circ}(\A) \quad = \quad
{\rm Ann}(V_i) \qquad \quad {\rm for} \,\, i =1,2,\ldots , r .$$
We conclude that $\, M(\B) +  {\rm Circ}(\A) \, = \,
{\rm Ann}\bigl( \{ V_1,\ldots,V_r \} \bigr)$,
which proves the identity (\ref{hypercogennn}).
\end{proof}

\begin{remark} We note that the above proof of
(\ref{crucialclaim}) is reversible, i.e.
Proposition~\ref{littlemma}
actually implies the levelness result of
Stanley \cite[Theorem III.3.4]{St1}
for matroids representable over $\Q$.
\end{remark}

We next rewrite the result of Theorem~\ref{hypercogenerator}
in terms of the other two presentations of our cohomology ring.
From  the perspective of the Lawrence toric variety $\,X(A^\pm, \theta) $,
it is most natural to work in a polynomial ring in $2n$ variables,
one for each torus-invariant divisor of $\,X(A^\pm, \theta)$.

\begin{corollary}
\label{2nvar}
The common cohomology ring  of the Lawrence toric variety and
the toric hyperk\"ahler variety has the presentation
$$
\,H^*(X(A^\pm,\theta);\Q) \quad = \quad
\Q[
\partial_{x_1},
\partial_{x_2},\ldots,
\partial_{x_n},
\partial_{y_1},
\partial_{y_2},\ldots,
\partial_{y_n} ] / {\rm Ann} \bigl(
 V_1(x-y),\ldots, V_r(x-y) \bigr). $$
\end{corollary}

\begin{proof}
The polynomials
$\, V_i(x-y) \, = \, V_i(x_1-y_1,x_2-y_2,\ldots,x_n-y_n)\,$
are annihilated precisely by the annihilators of $\,V_i(x)\,$
and by the extra ideal generators $\,
\partial_{x_1} +  \partial_{y_1} , \ldots,
\partial_{x_n} +  \partial_{y_n}$.
\end{proof}

This corollary states that the cogenerators of the Lawrence toric variety
are the volume polynomials of the maximal bounded faces of the associated
polyhedron $\,Q_\psi \, = \, P_\theta $. The same result holds
for any semi-projective toric variety, even if the maximal bounded faces
of its polyhedron have different dimensions. This can
be proved using Proposition~\ref{littlemma}.

The economical presentation of our cohomology ring is as a quotient
of a polynomial ring in $d$ variables $\partial_{t_1},\ldots,
\partial_{t_d}$. The matrix $A$ defines a surjective homomorphism
of polynomial rings
$$
 \alpha \,\, : \, \,
 \Q[\partial_{x_1},\ldots,\partial_{x_n}]
\, \,  \rightarrow \,\,
 \Q[\partial_{t_1},\ldots,\partial_{t_d}] \, , \quad
\partial_{x_j} \,\, \mapsto \,\,
\sum_{i=1}^d a_{ij} \partial_{t_i} , $$
and a dual injective  homomorphism
of polynomial rings
$$
 \alpha^* \,\, : \, \,  \Q[t_1,\ldots,t_d]  \,\, \rightarrow \,\,
 \Q[x_1,\ldots,x_n]  \,, \quad
t_i \,\, \mapsto \,\, \sum_{j=1}^n a_{ij} x_j . $$
The kernel of $\alpha$ equals $\,{\rm Circ}(\A)\,$ and
therefore
\begin{equation}
\label{tinycohomm}
 H^* \bigl( Y(A,\theta); \Q \bigr) \quad = \quad
\Q [ \partial_{t_1}, \ldots, \partial_{t_d} ] / \alpha ( M(\B)) .
\end{equation}
We obtain cogenerators for this presentation of
our cohomology ring as follows. Suppose that the indeterminate
vector $\, t = (t_1,\ldots,t_d) \,$ ranges over a small neighborhood
of $\,\theta = ( \theta_1,\ldots,\theta_d) \,$ in $\R^d$. For $t$
in this neighborhood, the polyhedron $P_t$ remains simple and
combinatorially isomorphic to $P_\theta$. The maximal bounded
faces of $P_t$ can be identified with $\Delta_1,\ldots,\Delta_r$
as before, but now the volume of $\Delta_i$ is a homogeneous
polynomial of degree $n-d$ in only $d$ variables:
$$ v_i(t) \quad = \quad v_i(t_1,\ldots,t_d) \quad = \quad
{\rm vol}(\Delta_i)
\quad {\rm for} \,\, i =1,2,\ldots , r .$$
The  polynomial $v_i(t)$ is the unique preimage of
the polynomial $V_i(x)$ under the
inclusion $\alpha^*$.

\begin{corollary}
\label{dvar}
The cohomology  of the Lawrence toric variety and
the toric hyperk\"ahler variety equals
$$ H^* \bigl(Y(A,\theta); \Q) \quad = \quad
\Q[ \partial_{t_1},\ldots,\partial_{t_d}]/
{\rm Ann} \bigl( \{v_1,\ldots,v_r\}). $$
\end{corollary}

\begin{proof}
A differential operator $\, f = f(\partial_{x_1},
\ldots,\partial_{x_n})$ annihilates $\,\alpha^*(v)\,$
for some $v = v(t_1,\ldots,t_d)$ if and only if
the operator $\alpha(f)$ annihilates $v$ itself.
This is the Chain Rule of Calculus. Hence
$$
 {\rm Ann}(\{v_1,\ldots,v_r\}) \,\, = \,\,
\alpha \bigl( {\rm Ann}(\{ V_1,\ldots,V_r \}) \bigr) \,\, = \,\,
\alpha \bigl( {\rm Circ}(\A) \, + \, M(\B) \bigr) \,\, = \,\,
\alpha(M(\B)). $$
The claim now follows from equation (\ref{tinycohomm}).
\end{proof}

\begin{remark} Since the cohomology ring of $Y(A,\theta)$ does not depend on
$\theta$, we get the remarkable fact that the vector space generated by the
volume polynomials does not depend on $\theta$ either.
\end{remark}

We close this section by presenting
an application  to combinatorics.  We use notation
and terminology as in \cite[Section III.3]{St1}.
Let $M $ be any matroid of
rank $n-d$ on $n$ elements which can be represented
over the field $\Q$, say, by a  configuration
$\B \subset \Z^{n-d}$ as above, and let
$ h(M) = (h_0,h_1,\ldots,h_{k})\,$ be its $h$-vector.
A longstanding open problem
is to characterize the $h$-vectors of matroids.
For a survey see \cite{CC} or \cite[Section III.3]{St1}. We wish to
argue that toric hyperk\"ahler geometry
can make a valuable contribution to this problem.
According to Corollary~\ref{hyperbetti}
the $h$-numbers of $M$ are precisely the
Betti numbers of the associated toric hyperk\"ahler
variety:
\begin{equation}
\label{hhhbetti}
 h_i(M) \quad = \quad {\rm rank} \, H^{2i}\bigl(Y(A,\theta);\Q \bigr) .
\end{equation}
As a first step, we prove the injectivity part of the
Hard Lefschetz Theorem for toric hyperk\"ahler varieties.
The $g$-vector of the matroid is  $\, g(M) =
(g_1,g_2,\dots,g_{{\lfloor {n-d\over 2}\rfloor}})$
where $g_i=h_{i}-h_{i-1}$.

\begin{theorem}
The $g$-vector of a rationally represented coloop-free matroid
is a Macaulay vector, i.e. there exists a graded
$\Q$-algebra $R = R_0 \oplus R_1 \oplus \cdots \oplus
R_{\lfloor {n-d\over 2}\rfloor}\,$ generated by $R_1$ and with
$g_i = {\rm dim}_\Q(R_i)$ for all $i$.
\label{application}
\end{theorem}

\begin{proof} Let $[D]\in H^2(Y(A,\theta);\Q)$ be the
class of an ample divisor. The restriction
 $D |_{X_j}$ to any component $X_j$ of the core is an
ample divisor on the projective toric orbifold $X_j$.
Consider the map
\begin{equation}
\label{mapL}
 L \,\, : \,\, H^{2i-2} \bigl( Y(A,\theta);\Q \bigr)
 \,\, \to \,\, H^{2i}\bigl(Y(A,\theta);\Q \bigr),
\end{equation}
given by multiplication with $[D]$.
We claim that this map is
 injective for $i=1,\dots,\lfloor {n-d\over 2}\rfloor $.
To see this, let $\alpha\in H^{2i-2}( Y(A,\theta);\Q)$ be a nonzero
cohomology class. Then according to equation
(\ref{crucialclaim}), there exists an index $j \in \{1,2,\ldots,r\}$
such that $\alpha | _{X_j}$ is
nonzero. Then the Hard Lefschetz Theorem for the projective
toric orbifold
$X_j$ implies that $\,\alpha | _{X_j} \cdot [D |_{X_j}] \,$
is a non-zero class in $\, H^{2i}( X_j;\Q)$. Its preimage
$\alpha \cdot [D]$ under the map $\phi_j$ is non-zero,
and we conclude that the map
(\ref{mapL}) is injective for $2i \leq n-d$.
Consider the quotient algebra
$\, R\, =\, H^*(Y(A,\theta);\Q)/\langle[D]\rangle$.
The injectivity result just established implies that
$$ g_i=h_{i}-h_{i-1}
\, = \, {\rm dim}_\Q
\bigl( H^{2i}(Y(A,\theta);\Q)/\langle[D]\rangle \bigr)
\, = \, {\rm dim}_\Q (R_i) . $$
 This completes the proof of  Theorem~\ref{application}.
\end{proof}

\begin{remark} After the submission of our paper we learned that Swartz
\cite{Sw} has given a different proof of Theorem~\ref{application} for all
coloop-free matroids.
\end{remark}

\section{Toric quiver varieties}
\label{quiver}

In this section we discuss an important class of toric hyperk\"ahler
manifolds, namely, Nakajima's quiver varieties in the
special case when the dimension vector has all coordinates
equal to one. Let $Q=(V,E)$ be a directed graph (a {\it quiver})
with $d+1$ vertices
$\,V=\{v_0,v_1,\dots,v_d\}$  and $n$ edges $\{e_{ij}:(i,j)\in E\}$.
We consider the group of all $\Z$-linear combinations
of $V$ whose coefficients sum to zero. We fix the basis
 $\{v_0-v_1,\dots,v_0-v_d\}$ for this group, which is
hence identified with $\Z^d$. We also identify $\Z^n$
with the group of $\Z$-linear combinations
$\sum_{ij} \lambda_{ij} e_{ij}$ of the set of edges $E$.
The boundary map of the quiver $Q$ is the following homomorphism
of abelian groups
\begin{equation}
\label{boundary}
 A \,\, : \,\, \Z^n \rightarrow \Z^d, \,\,
e_{ij} \, \mapsto \, v_i - v_j .
\end{equation}
Throughout this section we assume that
the underlying graph of $Q$ is connected. This
ensures that $A$ is an epimorphism. The kernel
of $A$ consists of all $\Z$-linear combinations
of $E$ which represent cycles in $Q$. We fix
an $n \times (n-d)$-matrix $B$ whose columns
form a basis for the cycle lattice ${\rm ker}(A)$.
Thus we are in the situation of (\ref{ses}).
The following result is well-known:

\begin{lemma}
The matrix $A$ representing the boundary map of a quiver $Q$
is unimodular.
\label{wellknown}
\end{lemma}

Every edge $e_{ij}$ of $Q$ determines one coordinate
function  $z_{ij}$ on $\C^n$
and two coordinate functions  $z_{ij},w_{ij}$ on $\H^n$.
The action of the $d$-torus on $\C^n$ and $\H^n$ given
by the matrix $A$ equals
\begin{equation}
\label{graphaction}
z_{ij} \mapsto t_i t_j^{-1} \cdot z_{ij}, \quad
w_{ij} \mapsto t_i^{-1} t_j \cdot w_{ij}.
\end{equation}
We are interested in the various quotients of
$\C^n$ and $\H^n$ by this action.
Since the matrix $A$ represents the quiver $Q$,
we  write $X(Q,\theta)$ instead of $X(A,\theta)$,
we write $X(Q^{\pm},\theta)$ instead of $X(A^{\pm},\theta)$,
and we write $\, Y(Q ,\theta)\,$ instead of  $\, Y(A ,\theta)$.
From Corollary 2.9 and Lemma~\ref{wellknown}, we conclude
that all of these quotients are manifolds when the
parameter vector $\theta$ is generic:

\begin{proposition}
Let $\theta $ be a generic vector in the lattice $\Z^d $.
Then $X(Q,\theta)$ is a smooth projective toric variety of dimension $n-d$,
$\, X(Q^{\pm} ,\theta)$ is a non-compact smooth toric variety
of dimension $2n-d$, and
$\, Y(Q ,\theta)$ is a smooth toric hyperk\"ahler variety
of dimension $2(n-d)$.
\end{proposition}

We call $Y(Q,\theta)$ a {\it toric quiver variety}. These
are precisely the quiver varieties of Nakajima \cite{Na}
in the case when the dimension vector has all coordinates equal to one.
Altmann and Hille \cite{AH} used the term ``toric quiver variety''
for the projective toric variety $X(Q,\theta)$, which is
a component in  the  common core of our
toric quiver variety $Y(Q,\theta)$ and its ambient
Lawrence toric variety $X(Q^{\pm}, \theta)$.
According to our general theory,
the manifolds $Y(Q,\theta)$ and  $X(Q^{\pm}, \theta)$
and their core
have the same integral cohomology ring, to be described
in terms of quiver data in Theorem~\ref{mainquiver}.

 Fix a vector $\theta \in \Z^d$ and a subset
$\tau \subseteq E$ which forms a {\em spanning tree} in $Q$.
Then there exists a unique linear combination
with integer coefficients $\,\lambda^\tau_{ij}\,$
which represents $\theta$ as follows:
$$ \theta \quad = \quad
\sum_{(i,j) \in \tau} \lambda^\tau_{ij} \cdot (v_i - v_j) . $$
Note that
the vector $\theta$ is generic if $\lambda_{ij}^\tau$
is non-zero for all spanning trees $\tau $
and all $(i,j) \in \tau $.

For every spanning tree $\tau$,
we define a subset of the monomials in $T = \C[z_{ij},w_{ij}]$
as follows.
$$
\sigma(\tau,\theta) \quad := \quad
\bigl\{ \,z_{ij} \, : \, (i,j) \in \tau \, \,
\hbox{and} \,\, \lambda_{ij}^\tau > 0 \, \bigr\} \,\, \cup \,\,
\bigl\{ \,w_{ij} \, : \, (i,j) \in \tau \, \,
\hbox{and} \,\, \lambda_{ij}^\tau < 0 \, \bigr\}
. $$
Recall that a {\it cut} of the quiver $Q$ is a collection $D$ of
edges which traverses a partition $\, (W, V \backslash W)\,$
of the vertex set $V$. We regard $D$ as a signed set by
recording the directions of its edges as follows
\begin{eqnarray*}
D^- \quad = \,\,\, &
\bigl\{ \, (i,j) \in E \, \, : \,\, i \in V \backslash W \,
\, \hbox{and} \,\,  j \in  W \, \bigr\}, \\
D^+ \quad = \,\,\, &
\bigl\{\, (i,j) \in E \, \, : \,\, i \in W \,\, \hbox{and} \,\,
  j \in V \backslash W \, \bigr\}.
\end{eqnarray*}
We now state our main result regarding toric quiver varieties:

\begin{theorem}
Let $\theta \in \Z^d$ be generic.
The Lawrence toric variety $X(Q^\pm, \theta)$ is the
smooth $(2n-d)$-dimensional toric variety
defined by the fan whose  $2n$ rays are the columns of
$\, \Lambda(\B)=\pmatrix{ {\bf I} & {\bf I} \cr
                {\bf 0} & B^T }$
and whose maximal cones are indexed by the sets
$\,\sigma(\tau,\theta)$, where $\tau$ runs
over all spanning trees of $Q$.
The toric quiver variety $Y(Q,\theta)$ is
the $2(n-d)$-dimensional submanifold of $X(Q^\pm, \theta)$ defined
 by the equations
$\, \sum_{(i,j) \in D^+} z_{ij} w_{ij}
=  \sum_{(i,j) \in D^-} z_{ij} w_{ij}\,$
where $D$ runs over all cuts of $Q$.
The common cohomology ring of these manifolds is the quotient
of $\,\Z[\partial_{ij} : (i,j) \in E]\,$ modulo the ideal
generated by the linear forms in $\, \partial \cdot B \,$ and
the monomials $\,\prod_{(i,j) \in D} \partial_{ij}\,$ where
$D$ runs over all cuts of $Q$.
\label{mainquiver}
\end{theorem}

A few comments are in place:
the variables $\,\partial_{ij}$, $ (i,j) \in E$,
are the coordinates of the row vector $\partial$,
so the entries of $\partial \cdot B$
are a cycle basis for $Q$. The equations
which cut out the toric quiver variety $\,Y(Q,\theta)\,$ lie in the
Cox homogeneous coordinate ring of
the Lawrence toric manifold $\,X(Q^{\pm},\theta)$.
 A more compact representation
is obtained if we replace ``cuts'' by
 ``cocircuits''.  By definition,
a {\it cocircuit} in $Q$ is a cut which is minimal
with respect to inclusion.
The proof of Theorem \ref{mainquiver} follows from our
general results for  integer matrices $A$
and is left to the reader.

Corollary \ref{hyperbetti} shows  that the
Betti numbers of $Y(Q,\theta)$ are the
$h$-numbers of the matroid of $B$.
This is not the usual graphic matroid of $Q$
but it is the {\em cographic matroid}
associated with $Q$.
Thus the Betti numbers of the   toric quiver variety $Y(Q,\theta)$
are the $h$-numbers of the cographic matroid of $Q$.
The generating function for the Betti numbers, the  $h$-polynomial
of the cographic matroid,
is known in combinatorics as the {\em reliability polynomial}
of the graph $Q$; see \cite{CC}.

\begin{corollary}
The Poincar\'e polynomial of the
toric quiver variety $Y(Q,\theta)$
equals the reliability polynomial of the graph $Q$,
which is the $h$-polynomial of its cographic matroid. In particular,
the Euler characteristic of  $Y(Q,\theta)$
coincides with the number of spanning trees of $Q$.
\end{corollary}

Recent work of Lopez \cite{ML}
gives an explicit enumerative interpretation of the
coefficients of the reliability polynomial of a graph and
hence of the Betti numbers of a toric quiver variety.
In particular, that paper proves
Stanley's longstanding conjecture on
$h$-vectors of matroid complexes
\cite[Conjecture III.3.6]{St1}
for the special case of cographic matroids.

\section{An example of a toric quiver variety}
\label{examples}

We shall describe a particular toric quiver variety
$\,Y(K_{2,3},\theta)\,$ of complex dimension four.
Consider the quiver in Figure 1,
the complete bipartite graph  $\,K_{2,3}\,$  given by
$\, d= 4 $, $\,  n=6 \,$ and $\, E = \{
(0,2),(0,3),(0,4),(1,2),(1,3),(1,4)\} $.

\setlength{\unitlength}{1.5cm}
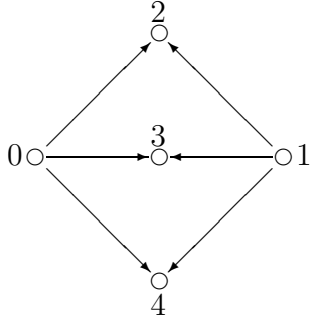
\begin{figure}[here]
\label{quiverQ}
\begin{center}
\begin{picture}(3,3)
\put(.4,1.5){\circle{0.15}}
\put(1.5,1.5){\circle{0.15}}
\put(2.6,1.5){\circle{0.15}}
\put(1.5,2.6){\circle{0.15}}
\put(1.5,0.4){\circle{0.15}}
\put(0.5,1.6){\vector(1,1){0.92}}
\put(0.5,1.5){\vector(1,0){0.90}}
\put(0.5,1.4){\vector(1,-1){0.92}}
\put(2.5,1.4){\vector(-1,-1){0.92}}
\put(2.5,1.5){\vector(-1,0){0.90}}
\put(2.5,1.6){\vector(-1,1){0.92}}
\put(0.15,1.43){0}
\put(1.42,1.6){3}
\put(2.71,1.43){1}
\put(1.42,2.7){2}
\put(1.42,.1){4}
\end{picture}
\end{center}
\caption{The quiver $K_{2,3}$}
\end{figure}
The matrix $A$ representing the boundary map (\ref{boundary})
is given in Figure 2.
The six columns of $A$ span the cone over a triangular prism
as depicted in Figure 3.
A Gale dual of this configuration is given by
the six vectors in the plane in  Figure 4.
The rows of $B^T$ span the cycle lattice of  $K_{2,3}$.

\begin{figure}[here]
\begin{center}
\begin{minipage}[b]{4.7cm}
\label{matrixA}
\begin{eqnarray*}A \, = \, \left[ \begin{array}{cccccc}
0 & 0 & 0 & -1 & -1 & -1 \\
1 & 0 & 0 & \phantom{-} 1 & \phantom{-}0 & \phantom{-}0 \\
0 & 1 & 0 & \phantom{-} 0 & \phantom{-} 1 & \phantom{-} 0 \\
0 & 0 & 1 & \phantom{-} 0 & \phantom{-} 0 & \phantom{-} 1 \\
 \end{array}\right]\end{eqnarray*}
\caption{The matrix A}
\end{minipage}
\hspace{1cm}
\begin{minipage}[b]{9cm}
\begin{center}
\label{configurationA}
\begin{picture}(6,4)
\put(1,2){\circle*{.15}}
\put(4,2){\circle*{.15}}
\put(2,3.5){\circle*{.15}}
\put(2,0.5){\circle*{.15}}
\put(5,3.5){\circle*{.15}}
\put(5,0.5){\circle*{.15}}
\put(1,2){\line(2,3){1}}
\put(1,2){\line(2,-3){1}}
\put(2,0.5){\line(0,1){3}}
\put(2,0.5){\line(1,0){3}}
\put(5,0.5){\line(0,1){3}}
 \put(2,3.5){\line(1,0){3}}
\multiput(1,2)(.8,0){4}{\line(1,0){.4}}
\multiput(4,2)(.4,.6){3}{\line(2,3){.2}}
\multiput(4,2)(.4,-.6){3}{\line(2,-3){.2}}
\put(0.2,1.8){$v_0-v_2$}
\put(1.7,3.6){$v_0-v_3$}
\put(1.7,.3){$v_0-v_4$}
\put(3.2,1.8){$v_1-v_2$}
\put(4.7,3.6){$v_1-v_3$}
\put(4.7,.3){$v_1-v_4$}
\end{picture}
\end{center}
\caption{The column vectors of the matrix $A$}
\end{minipage}
\end{center}
\end{figure}

\begin{figure}[here]
\begin{center}
\begin{minipage}[b]{6.5cm}
\label{matrixB}
\begin{eqnarray*} B^T \, =
\,\left[\begin{array}{cccccc} 1 &0 &-1 &-1& 0& 1 \\ 0&1&-1&0&-1&1\end{array}\right]
\end{eqnarray*}
\caption{Transpose of the matrix B}
\end{minipage}
\hspace{3cm}
\begin{minipage}[b]{6cm}
\begin{center}
\label{configurationB}
\pspicture(-2,-2)(2,2)
\psset{yunit=0.866026040cm}
\psline{->}(0,0)(2,0)
\psline{->}(0,0)(1,2)
\psline{->}(0,0)(-1,2)
\psline{->}(0,0)(-2,0)
\psline{->}(0,0)(-1,-2)
\psline{->}(0,0)(1,-2)
\rput(2.2,0){$b_{02}$}
\rput(1.2,2.2){$b_{14}$}
\rput(-.9,2.2){$b_{03}$}
\rput(-2.2,0){$b_{12}$}
\rput(-.9,-2.2){$b_{04}$}
\rput(1.2,-2.2){$b_{13}$}
\endpspicture
\end{center}
\caption{Rows of the matrix~$B$}
\end{minipage}
\end{center}
\end{figure}

Our manifolds are constructed algebraically from the
multigraded polynomial rings
\begin{equation}
\label{Sring}
 S \quad = \quad \C [ z_{02},z_{03},z_{04},z_{12},z_{13},z_{14} ]
\qquad \hbox{and} \quad
 T \quad = \quad S [ w_{02},w_{03},w_{04},w_{12},w_{13},w_{14} ],
\end{equation}
where the degrees of the twelve variables are given by the
columns of the matrix $A^{\pm}$:
\begin{equation}
\label{edgegrading}
{\rm degree}(z_{ij}) \quad = \quad - {\rm degree}(w_{ij})
\quad = \quad v_i - v_j .
\end{equation}
This grading corresponds to the torus action (\ref{graphaction})
on the polynomial rings $S$ and $T$.

Fix $\theta=(\theta_1,\theta_2,\theta_3,\theta_4)\in \Z^4$.
It represents the following linear combination of vertices
of $K_{2,3}$:
$$
(\theta_1 +
\theta_2 +
\theta_3 +
\theta_4 ) v_0
- \theta_1 v_1
- \theta_2 v_2
- \theta_3 v_3
- \theta_4 v_4
$$
The monomials $\,
z_{02}^{u_{02}}
z_{03}^{u_{03}}
z_{04}^{u_{04}}
z_{12}^{u_{12}}
z_{13}^{u_{13}}
z_{14}^{u_{14}}\,$
in the graded component $S_\theta$
correspond to the nonnegative $2 \times 3$-integer matrices
$\,\pmatrix{
u_{02} & u_{03} & u_{04} \cr
u_{12} & u_{13} & u_{14} \cr} \,$
with column sums $\,\theta_2,\theta_3,\theta_4\,$
and row sums $\, \theta_1 + \theta_2 + \theta_3 + \theta_4 \,$
and $\,-\theta_1 $.
For instance, for
$\,\theta = (-3,2,2,2)\,$
 there are precisely seven
monomials in $S_\theta$ as shown on Figure 6.
\begin{figure}[here]
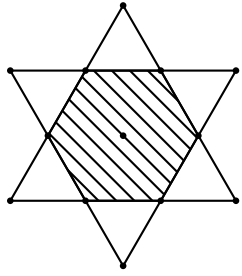

\label{theta}
\hspace{2cm}

\begin{minipage}[b]{6cm}
\hspace{2cm}
$$
\begin{array}{ll}
 S_\theta \, =\,
\C \, \bigl\{&
z_{02} z_{03} z_{04}
z_{12} z_{13} z_{14},\\ &
z_{02}^2      z_{04}
      z_{13}^2 z_{14},\\ &
z_{02}^2   z_{03}
      z_{13} z_{14}^2,\\ &
z_{02} z_{03}^2
z_{12}       z_{14}^2, \\ &
      z_{03}^2 z_{04}
z_{12}^2       z_{14}, \\ &
         z_{03} z_{04}^2
z_{12}^2 z_{13}, \\ &
z_{02}          z_{04}^2
z_{12} z_{13}^2 \quad \bigr\}
\end{array}
$$
\end{minipage}
\hspace{1.8cm}
\begin{minipage}[b]{4cm}
\begin{center}
\pspicture(0,0)(1,1)
\psset{yunit=0.866026040cm}
\pspolygon[fillstyle=vlines](1,0)(.5,1)(-.5,1)(-1,0)(-.5,-1)(.5,-1)
\pspolygon(0,2)(1.5,-1)(-1.5,-1)
\pspolygon(0,-2)(-1.5,1)(1.5,1)
\psdots*(1,0)(.5,1)(-.5,1)(-1,0)(-.5,-1)(.5,-1)(0,0)(0,2)(0,-2)(1.5,1)(-1.5,-1)(1.5,-1)(-1.5,1)
\endpspicture
\end{center}
\end{minipage}
\caption{Monomials in multidegree $\theta=(-3,2,2,2)$}
\end{figure}
Taking ``Proj'' of the algebra generated by these seven monomials
we get a smooth toric surface $X(K_{2,3},\theta)$ in $\P^6$. This surface
is the blow-up of $\P^2$ at three points.

As $\theta$ varies,  there are  eighteen different types of
smooth toric surfaces $X(K_{2,3},\theta)$.
They correspond to the eighteen chambers in the
triangular prism, or, equivalently, to the
eighteen complete fans on $\B$.
 This picture arises in the
 {\it Cremona transformation} of classical algebraic geometry,
where the projective plane is blown up at three points
and then the lines connecting them are blown down.
The eighteen surfaces are the intermediate blow-ups
and blow-downs.

We next describe the Lawrence toric varieties $\,X(K_{2,3}^{\pm},\theta)$
which are the GIT quotients of $\C^{12}$ by the action
(\ref{graphaction}). First, the (singular) affine quotient
$\,X(K_{2,3}^{\pm},0)  \,$ is the spectrum of the algebra
\begin{eqnarray*}
 T_0 \quad = \quad
& \C \bigl[
z_{02}  w_{02} ,
z_{03}  w_{03} ,
z_{04}  w_{04} ,
z_{12}  w_{12} ,
z_{13}  w_{13} ,
z_{14}  w_{14} , \,
 z_{02} z_{13} w_{12} w_{03},
\\
& z_{03} z_{12} w_{13} w_{02},
z_{02} z_{14} w_{12} w_{04},
z_{04} z_{12} w_{14} w_{02},
z_{03} z_{14} w_{13} w_{04},
z_{04} z_{13} w_{14} w_{03}
\bigr].
\end{eqnarray*}
This is the affine toric variety whose fan is the cone
over the $7$-dimensional Lawrence polytope given by the matrix
$\, \pmatrix{ {\bf I} & {\bf I} \cr
              {\bf 0} & B^T }$,
where ${\bf I}$ is the $ 6 \times 6$-unit matrix.
This Lawrence polytope has $160$ triangulations,
all of which are regular, so there are $160$ different
types of smooth Lawrence toric varieties $X(K_{2,3}^{\pm}, \theta)$
as $\theta$ ranges over the generic points in $\Z^4$.
For instance,  for $\theta = (-3,2,2,2)$ as in Figure 6,
$X(K_{2,3}^{\pm}, \theta)$ is constructed as follows.
The graded component $T_\theta$ is generated as a
$T_0$-module by $13$ monomials: the seven $z$-monomials
in $S_\theta$ and the six additional monomials:
\begin{equation}
\label{sixmonos}
w_{02} z_{03}^2 z_{04}^2 z_{12}^3, \,
w_{03} z_{02}^2 z_{04}^2 z_{13}^3, \,
w_{04} z_{02}^2 z_{03}^2 z_{14}^3, \,
w_{12} z_{13}^2 z_{14}^2 z_{02}^3, \,
w_{13} z_{12}^2 z_{14}^2 z_{03}^3, \,
w_{14} z_{12}^2 z_{13}^2 z_{04}^3.
\end{equation}
The $13$ monomial generators of $T_\theta$ correspond
to the $13$ lattice points in the star diagram
in Figure 6. The toric variety
$\, X(K_{2,3}^{\pm}, \theta) \, = \, {\rm Proj}( \oplus_{n \geq 0}
T_{n\theta})\,$ is characterized by its irrelevant ideal
in the Cox homogeneous coordinate ring $T$, which
is graded by (\ref{edgegrading}). The irrelevant ideal
is the radical of the monomial ideal $\langle T_\theta \rangle$.
It is generated by the $12$ square-free monomials obtained by
erasing exponents of the monomials in (\ref{sixmonos}) and Figure 6.
The $7$-simplices in the triangulation of the Lawrence polytope
are the complements of the supports of these twelve monomials,

We finally come to the toric quiver variety $Y(K_{2,3},\theta)$,
which is smooth and four-dimensional.
It is the complete intersection in the
Lawrence toric variety $\, X(K_{2,3}^{\pm}, \theta) \,$ defined by
the equations
$$
 z_{02} w_{02} +
 z_{03} w_{03} +
 z_{04} w_{04} \,\, = \,\,
z_{02} w_{02} + z_{12} w_{12} \,\, = \,\,
z_{03} w_{03} + z_{13} w_{13} \,\, = \,\,
z_{04} w_{04} + z_{14} w_{14} \,\, = \,\, 0 .$$
These equations are valid for all  $160$  toric quiver
varieties  $\,Y(K_{2,3}, \theta) $.
The cores of the manifolds vary  greatly.
For instance,
for $\theta = (-3,2,2,2)$, the core of $\, Y(K_{2,3}, \theta) \,$
consists of six copies of the projective plane $\P^2$ which are
glued to the blow-up of $\P^2$ at three points as in Figure 6.

The common cohomology ring of the $8$-dimensional Lawrence toric varieties
$\, X(K_{2,3}^{\pm}, \theta) \,$  and the $4$-dimensional toric quiver
varieties $\,  Y(K_{2,3}, \theta) \,$  is independent of $\theta$
and equals
\begin{eqnarray*}
\Z[\partial]/ \langle
 \partial_{03} \partial_{04} \partial_{12},
\partial_{02} \partial_{04} \partial_{13},
\partial_{02} \partial_{03} \partial_{14},
\partial_{13} \partial_{14} \partial_{02},
\partial_{12} \partial_{14} \partial_{03},
\partial_{12} \partial_{13} \partial_{04},
\\
 \partial_{02} \partial_{03} \partial_{04},
 \partial_{12} \partial_{13} \partial_{14},
 \partial_{02}  \partial_{12}, 
\partial_{03}  \partial_{13},
\partial_{04}  \partial_{14},
\,
\partial_{02} - \partial_{03} - \partial_{12} + \partial_{13} ,
\partial_{02} - \partial_{04} - \partial_{12} + \partial_{14}
\rangle.\end{eqnarray*}
From this presentation we can compute the Betti numbers as follows:
$$ H^*(Y(K_{2,3});\Z) \quad = \quad
 H^0(Y(K_{2,3});\Z) \, \oplus \,  H^2(Y(K_{2,3});\Z)
 \, \oplus \,  H^4(Y(K_{2,3});\Z)
 \quad = \quad
\Z^1 \, \oplus \, \Z^4 \, \oplus \, \Z^7. $$
The $7$-dimensional space of cogenerators is spanned
by the areas of the six triangles in Figure 6, e.g.,
$\, V_{\{03,04,12\}}(x) \,=\,
(x_{03} + x_{04} - x_{12})^2 $,
together with the  area polynomial of the hexagon
$$ V_{hex}(x) \,\, = \,\,
 2 x_{03} x_{14} +
 2 x_{14} x_{02} +
 2 x_{02} x_{13} +
 2 x_{13} x_{04}  +
 2 x_{04} x_{12} +
 2 x_{12} x_{03}
- x_{02}^2 - x_{03}^2 - x_{04}^2
- x_{12}^2 - x_{13}^2 - x_{14}^2. $$

\section{Which toric varieties  are hyperk\"ahler ?}
\label{ale}

Toric hyperk\"ahler varieties are constructed
algebraically as complete intersections in
Lawrence toric varieties, but they are generally
not toric varieties themselves. What we mean by this is that
there does not exist a subtorus of the
dense torus of $X(A^\pm,\theta)$ such that
$Y(A,\theta)$ is an orbit closure of that subtorus.
The objective of this section is to characterize
and study the rare exceptional cases when $Y(A,\theta)$
happens to be a toric variety. We are particularly
interested in the case of manifolds, when $A$ is
unimodular. The following is the  main result in this section.

\begin{theorem}
\label{ALEthm}
A toric manifold is a toric hyperk\"ahler variety
if and only if it is a  product of
ALE spaces of type $A_n$ if and only if
it is a toric  quiver variety $X(Q,\theta)$
where $Q$ is a disjoint union of cycles.
\end{theorem}

The ALE space of type $A_n$ is denoted
$\,\C^2 / \! / \Gamma_n \,$ where $\Gamma_n$ is the cyclic
group of order $n$ acting on $\C^2$ as
the matrix group $\, \bigl\{ \pmatrix{ \eta & 0 \cr 0 & \eta^{-1}} \, : \,
 \, \eta^n =  1 \, \bigr\} $. The smooth surface
$\,\C^2 / \! / \Gamma_n \,$
 is defined as the unique crepant resolution
of the $2$-dimensional cyclic quotient singularity
$$ \C^2 / \Gamma_n \quad = \quad
{\rm Spec}\, \C[x,y]^{\Gamma_n} \quad = \quad
{\rm Spec}\, \C[x^n,xy, y^n ] . $$
Equivalently, we can construct
 $\,\C^2 / \! / \Gamma_n \,$  as the
smooth toric surface whose fan $\Sigma_n$
consists of the cones
$\, \R_{\geq 0} \{ (1,i-1), \, (1,i) \}\,$
for $i = 1,2,\ldots,n\,$ and whose lattice is
the standard lattice $\Z^2$.

Let us start out by showing that the ALE space
$ \C^2 /\!/ \Gamma_n $ is indeed a toric quiver variety.
Let $C_n$ denote the {\it $n$-cycle}. This is the quiver with
vertices $V = \{0,1,\ldots,n-1 \}$ and edges
$$ E \,\, = \,\, \bigl\{ \,
(0,1), \, (1,2), \, (2,3), \,
\ldots, (n \! - \! 2, n-1), \, (n-1, 0) \,\bigr\}. $$
We prove the following well-know result to illustrate our constructions
in this paper.

\begin{lemma}
\label{cyclelem}
The affine quiver variety $Y(C_n,0)$ is isomorphic to
$ \C^2 / \Gamma_n \,$ and for any generic
vector $\theta \in \Z^{n-1}$,
the  smooth  quiver variety $Y(C_n,\theta)$
is isomorphic to the ALE space
$\, \C^2 /\!/ \Gamma_n $.
\end{lemma}

\begin{proof}
The boundary map of the $n$-cycle
 $C_n$ has the format $\, \Z^n \rightarrow \Z^{n-1}$
and looks like
$$ A \quad = \quad
\pmatrix{
 1 & -1 & 0  & 0 & \cdots & 0  \cr
 0 &  1 & -1 & 0 & \cdots & 0 \cr
 0 &  0 &  1 & -1 &  \cdots &  0 \cr
 \vdots &  \vdots & \ddots &   \ddots &  \ddots & \vdots \cr
 0 &  0 &  0 &    \cdots & 1 & -1 \cr}. $$
and its Gale dual is the $1 \times n$-matrix with all entries
equal to one:
\begin{equation}
\label{Bisallones}
 B^T \quad  = \quad \pmatrix{ 1 &  1 &  1 & \cdots & 1 \cr}.
\end{equation}
The torus $\T^{n-1}_\C$ acts via $A^\pm$ on the
 polynomial ring
$\, T \, = \, \C \, [\, z_{i,i+1}, w_{i,i+1} \, : \,
i = 0,\ldots,n-1 \,] $.
The affine Lawrence toric variety $\,X(C_n^\pm,0) =
\C^{2n} /\!/_0 \T^{n-1}_\C \,$ is the
spectrum of the invariant ring
$$ T_0 \quad = \quad \C \,\bigl[ \,
 z_{01}  w_{01} \,, \,
 z_{12}  w_{12} \,, \, \ldots, \,
 z_{n-1,0}  w_{n-1,0} \,, \,\,
z_{01} z_{12} \cdots z_{n-1,0} \,, \,\,
w_{01} w_{12} \cdots w_{n-1,0} \,\bigr]. $$
The common defining ideal  of all the quiver
varieties $Y(C_n,\theta)$ is the following ideal in $T$:
$$ {\rm Circ}(\B) \quad = \quad
\langle \,
z_{i-1,i} w_{i-1,i} \,- \,
z_{i,i+1} w_{i,i+1} \, \, : \,\,
i = 1,2,\ldots,n \, \rangle . $$
All indices are considered modulo $n$.
The quiver variety
$Y(C_n,0)$ is the spectrum of
$\,T_0 / (T_0 \,\cap \, {\rm Circ}(\B))$.
Dividing $T_0$ by
$\,T_0 \,\cap \, {\rm Circ}(\B)\,$ means
erasing the double indices of all variables:
$$ \,T_0 / (T_0 \,\cap \, {\rm Circ}(\B)) \quad \simeq \quad
\C[\, z w , \,z^n, \, w^n \, ] . $$
Passing to the spectra of these rings
proves our first assertion:
$\,Y(C_n,0) \simeq \C^2/\Gamma_n$.

For the second assertion, we first note that
$\, \theta = (\theta_1,\ldots,\theta_{n-1})\,$
is generic for $\A^\pm$ if and only if all
consecutive coordinate sums
$\,  \theta_i +
\theta_{i+1} +  \cdots + \theta_{j} \,$ are non-zero.
The associated hyperplane arrangement $\Gamma(\A)$
is linearly isomorphic to the braid arrangement
$\, \{ u_i = u_j \}$. It has $n\,!\,$ chambers,
and the symmetric group acts transitively on
the chambers.
Hence it suffices to prove our claim
$\,Y(C_n,\theta)\,\simeq \, \C^2/\!/ \Gamma_n\,$
for only one vector $\theta$
which lies  in the interior of any chamber.

We fix the generic vector $\, \theta = (1,1,\ldots,1)$.
There are  $n$ monomials of degree $\theta$ in $T$, namely,
\begin{equation}
\label{nmonos}
\prod_{j=1}^{i-1} z_{j-1,j}^{i-j}
\cdot
\prod_{k=i+1}^n w_{k-1,k}^{k-i} \qquad
\qquad \hbox{for} \,\,\, i = 1,2,\ldots,n .
\end{equation}
The images of these monomials are minimal generators
of the  $\,T_0 / (T_0 \,\cap \, {\rm Circ}(\B))$-algebra
$$ \bigoplus_{r=0}^\infty
T_{r \theta} / (T_{r \theta} \,\cap \, {\rm \Circ}(\B)) .$$
By definition, $Y(C_n,\theta)$ is the projective
spectrum of this $\N$-graded algebra. Applying our isomorphism
``erasing double indices'', the images of our $n$ monomials
in (\ref{nmonos}) translate into
\begin{equation}
\label{nbivariatemonos}
z^{i\choose  2} \cdot
w^{n-i+1 \choose 2}  \qquad
\qquad \hbox{for} \,\,\, i = 1,2,\ldots,n .
\end{equation}
Hence $Y(C_n,\theta)$ is the projective spectrum
of the $\C[zw,z^n,w^n]$-algebra generated
by (\ref{nbivariatemonos}). It is straightforward
to see that this is the toric surface
with fan $\Sigma_n$, i.e.~the
ALE space $\C^2/\!/ \Gamma_n$.
\end{proof}

It is instructive to write down our presentations
for the  cohomology ring of the ALE space
$\,Y(C_n,\theta) = \C^2 / \! / \Gamma_n $.
The circuit ideal  of the $n$-cycle  is the principal ideal
$$ {\rm Circ}(\A) \quad = \quad
\langle\,  \partial_{01} +   \partial_{12} +
 \partial_{23} +  \cdots + \partial_{n-1,0} \, \rangle.  $$
The matroid ideal $M(\B) $ is
generated by all quadratic squarefree monomials
in $\Z[\partial ]$.
It follows that
$\,\Z[ \partial]/({\rm Circ}(\A) + M(\B))\,$ is isomorphic to
a polynomial ring in $n-1$ variables modulo the square
of the maximal ideal generated by the variables, and hence
$$ H^*(Y(C_n,\theta); \Z) \quad = \quad
H^0(Y(C_n,\theta); \Z) \,\,\, \oplus \,\,\,
H^2(Y(C_n,\theta); \Z)
\quad \simeq \quad \Z^1 \, \oplus \, \Z^{n-1}. $$

On our way towards proving Theorem~\ref{ALEthm},
let us now fix an epimorphism
$\,A : \Z^n \rightarrow \Z^d \,$ and a generic vector
$\,\theta \in \Z^d $.  We assume that
 $\A$ is not a cone, i.e.
the zero vector is not in $\B$.
We do not assume  that $A$ is unimodular.
By a {\it binomial} we mean a polynomial with two terms.

\begin{proposition}
\label{propi898}
The following three statements are equivalent:
\begin{itemize}
\item[(a)] The hyperk\"ahler toric variety $\,Y(A,\theta)\,$
is a toric subvariety of $X(A^\pm,\theta)$.
\item[(b)] The ideal $\, {\rm Circ}(\B)\,$ is generated by binomials.
\item[(c)] The configuration $\B$ lies on
$n-d$ linearly independent lines through the origin in $\R^{n-d}$.
\end{itemize}
\end{proposition}

\begin{proof}
The condition (b) holds if and only if
the matrix $A$ can be chosen to have
two nonzero entries in each row. This defines
a graph ${\cal G}$ on $\{1,2,\ldots,n\}$, namely,
$j$ and $k$ are connected by an edge
if there exists $i \in \{1,\ldots,d\}$ such that
 $\, a_{ij} \not= 0 \,$ and
 $\, a_{ik} \not= 0 $.
The graph ${\cal G}$ is a disjoint union of $n-d$ trees.
Two indices  $j$ and $k$ lie in the same connected component
of ${\cal G}$ if and only if the vectors $b_j$ and $b_k$
are linearly dependent. This shows that (b) is equivalent to (c).

Suppose that (b) holds. Then
the prime ideal ${\rm Circ}(\B)$ is generated
by the quadratic binomials
$\,a_{ij} z_j w_j + a_{ik} z_k w_k\, $ indexed
by the edges $(j,k)$ of ${\cal G}$.
The corresponding coefficient-free equations
$$ z_j w_j \,\,\, = \,\,\, z_k w_k \qquad
{\rm for} \,\,\, (j,k) \in {\cal G}. $$
define a subtorus $\,\T\,$ of the dense torus
of the Lawrence toric variety $X(A^\pm,\theta)$, and
the equations
$$ a_{ij} z_j w_j \,\, + \,\, a_{ik} z_k w_k \quad  = \quad 0
\qquad {\rm for} \,\,\, (j,k) \in {\cal G}. $$
define an orbit of $\T$ in the dense torus of  $X(A^\pm,\theta)$.
The solution set of the same equations in $X(A^\pm,\theta)$
has the closure of that $\T$-orbit as  one of its irreducible
components. But that solution set is our hyperk\"ahler variety
 $Y(A,\theta)$. Since $Y(A,\theta)$ is irreducible, we can conclude
that it coincides with the closure of the $\T$-orbit. Hence
$Y(A,\theta)$ is a toric variety, i.e. (a) holds.

For the converse, suppose that (a) holds.
The irreducible subvariety $Y(A,\theta)$ is defined
by a homogeneous prime  ideal $J$ in
the homogeneous coordinate ring $T$ of $X(A^\pm,\theta)$.
Since $Y(A,\theta)$ is a torus orbit closure,
the ideal $J$ is generated by binomials.
The ideal  ${\rm Circ}(\B)$ has the same zero set as
$J$ does, and therefore, by the Nullstellensatz and results of Cox,
$\, {\rm rad}\bigl({\rm Circ}(\B) : B_\theta^\infty \bigr) \, = \, J$.
Our hypothesis $0 \not\in \B$ ensures that
${\rm Circ}(\B)$ itself is a prime ideal, and therefore
we conclude $\, {\rm Circ}(\B) \, = \, J$. In particular,
this ideal is generated by binomials, i.e. (b) holds.
\end{proof}

\smallskip

\noindent {\sl Proof of Theorem~\ref{ALEthm}}:
Suppose that $Q$ is a  quiver with
connected components $Q_1,\ldots,Q_r$.
Then its boundary map is given by
a matrix with block decomposition
\begin{equation}
\label{Adecomp}
 A \quad = \quad A_1 \, \oplus \,
A_2 \, \oplus \, \cdots \,\oplus \, A_r ,
\end{equation}
where $A_i$ is the boundary map of $Q_i$.
There is a corresponding  decomposition of
the Gale dual
\begin{equation}
\label{Bdecomp}
 B \quad = \quad B_1 \, \oplus \,
B_2 \, \oplus \, \cdots \,\oplus \, B_r .
\end{equation}
In this situation, the toric hyperk\"ahler
variety $Y(A,\theta)$ is the direct product
of the toric hyperk\"ahler varieties
$Y(A_i, \theta)$ for $i=1,\ldots,r$.
For our quiver $Q$ this means
$$
 Y(Q,\theta) \quad = \quad
 Y(Q_1,\theta) \, \times \,
 Y(Q_2,\theta) \, \times \,
\cdots \, \times \,   Y(Q_r,\theta). $$
Using Lemma \ref{cyclelem}, we conclude that
a manifold is   a product of
ALE spaces of type $A_n$ if and only if
it is a toric  quiver variety $Y(Q,\theta)$
where $Q$ is a disjoint union of cycles $\, C_{n_i}$.

The matrix $A$ in (\ref{Adecomp}) is unimodular if and only if
the matrices $A_1,\ldots,A_r$ are unimodular.
Hence a  product of
toric hyperk\"ahler  manifolds
is  a toric hyperk\"ahler manifold.
In particular, a  product of ALE spaces
$\, \C^2/\Gamma_{n_i} \,$ is a toric hyperk\"ahler
manifold which is also a toric variety.

For the converse, suppose that
$Y(A,\theta)$ is a toric hyperk\"ahler
manifold which is also a toric variety, so that
statement (a) in Proposition \ref{propi898} holds.
Statement (c) in Proposition \ref{propi898} says that
the matrix $B$ has a decomposition  (\ref{Bdecomp})
where $r = n-d$ and each $B_i$ is a matrix
with exactly one column. We may assume that none
of the entries in $B_i$ is zero.  The Gale dual $A_i$
of $B_i$ is a unimodular matrix, and hence $B_i$ is unimodular.
For a matrix with one column this means that all entries in
$B_i$ are either $+1$ or $-1$. After trivial sign changes,
this means
$\, B_i^T \, = \, ( \, 1 \, 1\, \ldots \, 1 \,)$.
Now we are in the situation of (\ref{Bisallones}),
which means that $\,Y(A_i,\theta)\,$ is an ALE space
$\, \C^2/\Gamma_{n_i} $. $\square$

\end{document}